\documentclass{mcom-l}
\usepackage{amssymb}
\usepackage{graphicx}
\usepackage{mathtools}
\usepackage{color}
\usepackage{cite}

\copyrightinfo{2014}{American Mathematical Society}
\newtheorem{theorem}{Theorem}[section]
\newtheorem{lemma}[theorem]{Lemma}
\newtheorem{coro}[theorem]{Corollary}
\newtheorem{proposition}[theorem]{Proposition}
\usepackage[normalem]{ulem}

\theoremstyle{definition}
\newtheorem{definition}[theorem]{Definition}

\newtheorem{defn}{Definition}

\theoremstyle{remark}

\newtheorem{rem}{Remark}

\numberwithin{equation}{section}
\newcommand{\ZZ}{\mathbb{Z}}
\newcommand{\RR}{\mathbb{R}}
\newcommand{\NN}{\mathbb{N}}
\newcommand{\CC}{\mathbb{C}}

\newcommand{\Bt}{{\tilde B}}
\newcommand\ba{\mathbf{a}}

\newcommand\bq{\mathbf{q}}
\newcommand\bv{\mathbf{v}}
\newcommand\bw{\mathbf{w}}
\newcommand\bbf{\mathbf{f}}

\def\bcalA{\mbox{\boldmath ${\mathcal A}$}}

\newcommand\bM{\mathbf{M}}

\newcommand\ra{{\rm a}}

\newcommand\ri{{\rm i}}

\newcommand\ie{{\it\thinspace i.e.}\ }

\begin{document}

\title[Exponential splines and pseudo-splines]{Exponential splines and pseudo-splines: generation versus reproduction of\\ exponential polynomials}

\author[C. Conti]{Costanza Conti}
\address{Dipartimento di Ingegneria Industriale, Universit\`{a} di Firenze,
Viale Morgagni 40/44, 50134 Firenze, Italy}
\email{costanza.conti@unifi.it}

\author[L. Gemignani]{Luca Gemignani}
\address{Dipartimento di Informatica,
Universit\`{a} di Pisa, Largo Bruno Pontecorvo, 3
56127 Pisa, Italy}
\email{l.gemignani@di.unipi.it}

\author[L.Romani]{Lucia Romani}
\address{Dipartimento di Matematica e Applicazioni, Universit\`{a} di Milano-Bicocca, Via R. Cozzi 55, 20125 Milano, Italy}
\email{lucia.romani@unimib.it}


\date{}

\dedicatory{}

\begin{abstract}
Subdivision schemes are iterative methods for the design of smooth curves and surfaces.
Any linear subdivision scheme can be identified by a sequence of Laurent polynomials, also called subdivision symbols,
which describe the linear rules determining successive refinements of coarse initial meshes.
One important property of subdivision schemes is their capability of exactly reproducing in the limit  specific types of functions from which the data is sampled. Indeed, this property is linked to the approximation order of the scheme and to its regularity. When the capability of reproducing polynomials is required, it is possible to define a family of subdivision schemes that allows to meet various demands for balancing approximation order, regularity and support size. The members of this family are known in the literature with the name of pseudo-splines. In case reproduction of exponential polynomials instead of polynomials is requested, the resulting family turns out to be
 the non-stationary counterpart of the one of pseudo-splines, that we here call the family of exponential pseudo-splines.
The goal of this work is to derive the explicit expressions of the subdivision symbols of
exponential pseudo-splines and to study their symmetry properties as well as their convergence and regularity.
\end{abstract}

\maketitle

\section{Introduction}
\begin{figure}[h!]
\centering
\hspace{-0.5cm}
{\includegraphics[trim= 5mm 5mm 5mm 5mm, clip, width=4.0cm]{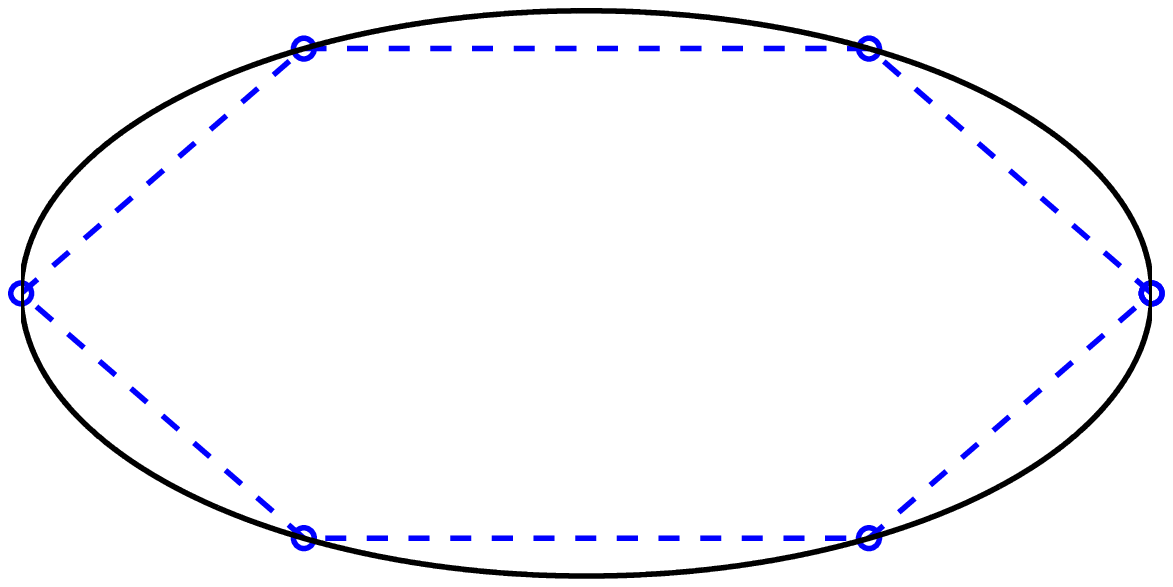}}\hspace{0.1cm}
{\includegraphics[trim= 5mm 5mm 5mm 5mm, clip, width=4.0cm]{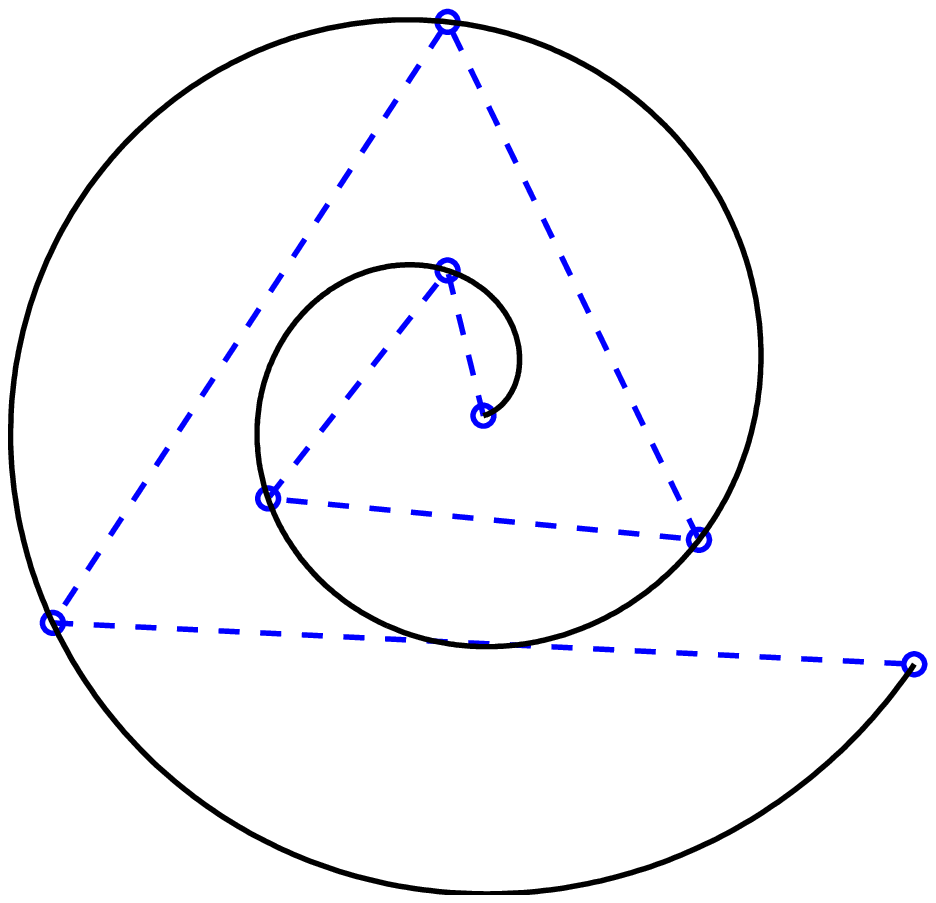}}\hspace{-0.1cm}
{\includegraphics[trim= 5mm 5mm 5mm 5mm, clip, width=4.0cm]{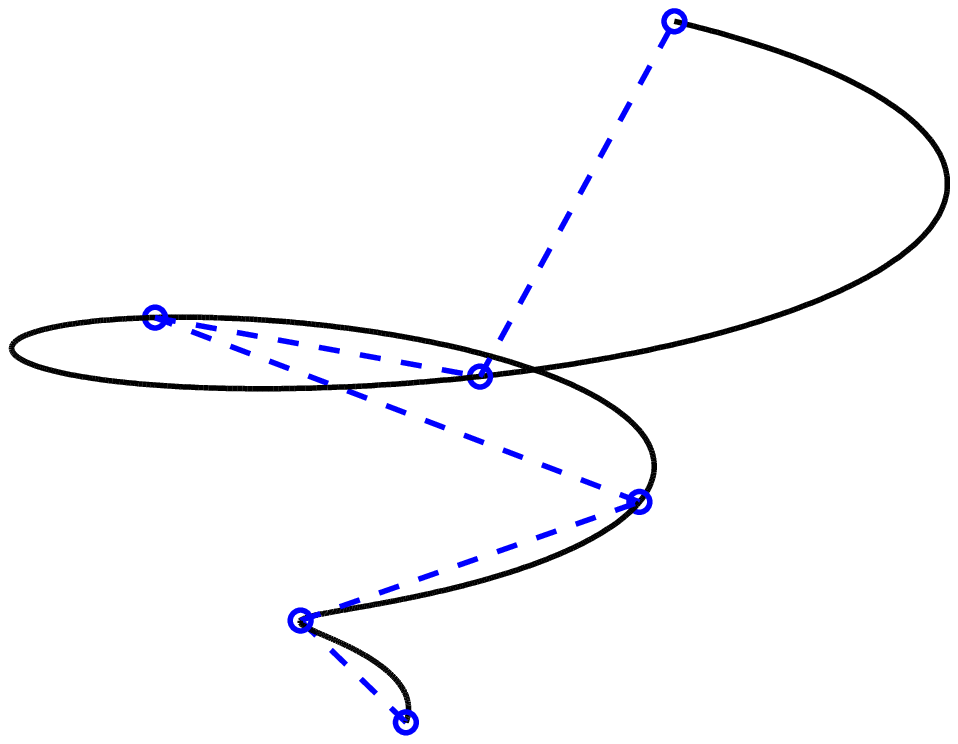}}
\caption{
Examples of reproduction of shapes
described by exponential polynomials using non-stationary subdivision schemes: initial control points and corresponding limit shapes for the subdivision schemes in \cite{ContiRomani11}. }
\label{figA}
\end{figure}

Subdivision schemes are efficient tools for the design of smooth curves and surfaces in many applicative areas such as
computer--aided geometric design, curve and surface reconstruction, signal/image processing.
Since in all these areas the capability of representing shapes described by polynomial, trigonometric or hyperbolic functions is fundamental (see Figure \ref{figA}), interpolating and approximating subdivision schemes based on exponential B--splines and inheriting their generation properties, have been recently introduced \cite{BCR07, BCR10, CC13, CCR13, CGR09, CGR11, ContiRomani09, ContiRomani11, DLL03, R09}.
The property of reproduction of exponential polynomials is also important since strictly connected to the approximation order of subdivision schemes and to their regularity \cite{ContiRomaniYoon2014}. In fact, the higher is the number of exponential polynomials reproduced, the higher is the approximation order and the possible regularity of the scheme.  Indeed, in application, we aim at subdivision schemes with exponential polynomial reproduction properties, that allow to meet various demands for balancing approximation order, regularity and support size. Such kind of schemes turn out to constitute  the family of \emph{exponential pseudo-splines}, the non-stationary counterpart of polynomial pseudo-splines. The latter family neatly fills the gap between B-splines and $2n$-point interpolatory subdivision schemes -both extreme cases of pseudo-splines: while B-splines stand out due to their high smoothness and short support, they provide a rather poor approximation order; in contrast, the limit functions of $2n$-point interpolatory subdivision schemes have optimal approximation order but low smoothness and large support.

Binary, primal and dual (polynomial) pseudo-splines -the first originally presented by Dong and Shen \cite{DongShen07}, the latter successively discovered by Dyn et al. \cite{DHSS08} and generalized to
any arity and to arbitrary parametrizations by Conti and Hormann \cite{ContiHormann11}- are both obtained by means of stationary subdivision schemes whose symbols can be read as a suitable
polynomial ``correction'' of the order-$N$ polynomial B-spline symbol $B_{N}(z)$, since of the form $a_{M,N}(z)=B_{N}(z)c_M(z)$. The polynomial correction $c_M(z)$ is such that the subdivision schemes with symbols $a_{M,N}(z)$ are the ones of minimal support that, besides \emph{generating polynomials of degree $N-1$}, satisfy the conditions for \emph{reproduction of polynomials of degree $M-1$}, with $M\le N$. We recall that while with generation we mean the subdivision capability to provide specific type of limit functions, with reproduction we mean  the capability of a subdivision scheme to reproduce in the limit exactly the same function from which the data is sampled.
Similarly to the stationary case, we here define exponential pseudo-splines by means of $k$-level subdivision symbols which are a suitable ``correction'' of the $k$-level subdivision symbols $B^{(k)}_{N,\Gamma}(z)$ of exponential B-splines, i.e. of the form $a^{(k)}_{M,N,\Gamma}(z)=B^{(k)}_{N,\Gamma}(z) c_{M,\Gamma}^{(k)}(z)$. Here $\Gamma$ identifies the particular space of exponential polynomials $EP_{\Gamma}$ we deal with, while $N$ and $M$ are related to the number of exponential polynomials that are being generated and reproduced, respectively. Again, $c_{M,\Gamma}^{(k)}(z)$ is such that
the symbols $a_{M,N,\Gamma}^{(k)}(z)$ are of minimal support and satisfy the conditions for reproduction of the space of exponential polynomials $EP_{\Gamma}$ generated by the exponential B-splines with symbols $B^{(k)}_{N,\Gamma}(z)$, or a subset of it.

\smallskip \noindent
The main contribution of this paper consists in showing how the symbols of exponential pseudo-splines can be explicitly derived.
Indeed, we provide the expressions of the inverse matrices of the linear systems arising by imposing the algebraic conditions for exponential polynomial reproduction which were firstly given in \cite{ContiRomani11} and successively extended to any arbitrary arity in \cite{CCR13}. We also show that, under the
\emph{symmetry assumption on $\Gamma$} (or on its subset), the symbol $a_{M,N,\Gamma}^{(k)}(z)$ has the same symmetry as $B_{N,\Gamma}^{(k)}(z)$. To prove the latter we also discover remarkable algebraic properties, never highlighted so far, of symmetric non-stationary subdivision symbols. As a minor contribution, we show how the $k$-level normalization factor of the exponential B-spline symbol
can be selected in accordance with the shift parameter in order to ensure that the exponential B-spline is \emph{correctly normalized}, namely,  besides generating the space $EP_{\Gamma}$, it \emph{reproduces} a specific pair of exponential polynomials $\{e^{\theta_{\ell} \, x},\ e^{-\theta_{\ell} \, x}\} \in EP_{\Gamma}$.
Finally, we additionally provide a convergence and regularity analysis of the non-stationary subdivision schemes corresponding to the exponential pseudo-spline symbols here derived. This is possible by first showing that exponential pseudo-splines are asymptotically similar to polynomial pseudo-splines, and then combining recent advances on convergence and regularity of non-stationary subdivision schemes presented in \cite{CharinaContiGuglielmiProtasov2014} and in \cite{ContiRomaniYoon2014}.

\medskip \noindent
The remainder of the paper is organized as follows. In Section \ref{section_base} we recall basic notions on non-stationary subdivision schemes reproducing exponential polynomials. Then, in Section \ref{sec:symmetry} we discuss new important results concerning symmetry properties of such subdivision schemes. Symmetric exponential B-splines are recalled in Section \ref{exponential B-splines} where
accordance between their parameter shift and their normalization factor is also considered with respect to their reproduction capabilities.
The derivation of the symbols of exponential pseudo-splines is provided in Section \ref{L-exponential pseudo-splines} where the symmetry properties of such symbols are also discussed. Convergence and regularity of the new family of (non-stationary) exponential pseudo-spline subdivision schemes are then investigated in Section \ref{convergence}. As an example of application of the presented theoretical results, the expression of the subdivision symbols of a new family of exponential pseudo-splines is also explicitly derived in Section \ref{Example}, where pictures of basic limit functions of the corresponding subdivision schemes are also given. The closing Section \ref{conclusion} is to draw conclusions.

\section{Non-stationary subdivision schemes and exponential polynomials }\label{section_base}
This paper deals with non-stationary subdivision schemes and reproduction of exponential polynomials.
The interest in non-stationary subdivision schemes arose in the last ten years after it
was pointed out that they are able to reproduce conic sections, spirals or widely used trigonometric/hyperbolic
curves and surfaces, as well as they are featured by tension parameters that allow, on
the one side, to obtain considerable variations of shape and, on the other side, to get
close to the initial mesh as much as desired (see \cite{ContiRomani09, Warren2002, R09}). Since numerical methods based on subdivision schemes are relatively simple to implement
and highly intuitive in use, they are currently widely exploited in modeling freeform curves and surfaces in computer games and animated movies.
The potential of subdivision schemes  has recently become apparent also in the context of
Isogeometric Analysis (IgA), a modern computational approach that offers the possibility
of integrating finite element analysis (FEA) into conventional CAD systems (see, e.g.,
\cite{Burkhart10,COS00,CSAOS02}).
However, the employment of IgA in conjunction with subdivision schemes is nowadays only
restricted to the class of stationary methods.
This is due to the fact that non-stationary subdivision schemes still require the development of further theoretical results that turn out to be fundamental to support their practical use.\\

\smallskip \noindent
Following the notation in \cite{DL02}, for any $k\ge 0$ we denote by $\ba^{(k)}:=\{\ra^{(k)}_j,\ j\in \ZZ\}$ the finite set of real coefficients corresponding to the so called {\sl $k$-level mask} of a non-stationary subdivision scheme and we define by $
a^{(k)}(z):=\sum_{j \in \ZZ}, \ra^{(k)}_j z^j,  \ z \in \CC \backslash\{0\}\,,
$
the Laurent polynomial whose coefficients are exactly the entries of $\ba^{(k)}$.
The previous polynomial is commonly known as the {\sl k-level symbol}
of the non-stationary subdivision scheme. With any mask $\ba^{(k)}$ comes a linear subdivision operator identifying a refinement process, that is the process which transforms a set of real data at level $k$, $\bbf^{(k)} =\{f^{(k)}_i \in \RR,\ i\in \ZZ\}$, into the denser set $\bbf^{(k+1)}$ given by
\begin{equation}\label{def:suboper}
\bbf^{(k+1)}:=S_{\ba^{(k)}}\bbf^{(k)}, \quad \hbox{where}\quad (S_{\ba^{(k)}}\bbf^{(k)})_i:=\sum_{j\in \ZZ} \ra^{(k)}_{i-2j} \, f_j^{(k)}, \quad \ \forall \, k\ge 0.
\end{equation}

\noindent The \emph{subdivision scheme} consists in the repeated application of the subdivision operators starting from any initial ``data" sequence $\bbf^{(0)} \equiv \bbf:=\{f_i \in \RR,\ i\in \ZZ\}$, and therefore is shortly denoted by $\{S_{\ba^{(k)}},\ k\ge 0\}$.

\smallskip \noindent
Since  the subdivision process generates denser and denser sequences of data, attaching the data $f_i^{(k)}$ generated at the $k$-th step to the \emph{parameter} values $t^{(k)}_i$ with
$t^{(k)}_i<t_{i+1}^{(k)}$ and $t_{i+1}^{(k)}-t_i^{(k)}=2^{-k},\ k\ge 0$, a notion of convergence can be established by
taking into account the piecewise linear function $F^{(k)}$ that interpolates the data (namely $
  F^{(k)}(t_i^{(k)}) = f_i^{(k)}, \
  F^{(k)}|_{[t_i^{(k)},t_{i+1}^{(k)}]} \in \Pi_1, \
  i\in\ZZ, \  k\geq0
$).
If
the sequence of continuous functions $\{F^{(k)},\ k\ge 0\}$ converges uniformly, we denote its
limit by
 $$ g_\bbf := \displaystyle{\lim_{k\to +\infty} S_{\ba^{(k)}}S_{\ba^{(k-1)}} \cdots S_{\ba^{(0)}} {\bf f}=\lim_{k\to +\infty} F^{(k)}},$$
and say that $g_\bbf$ is the \emph{limit function} of
the non-stationary subdivision scheme based on the rules in (\ref{def:suboper}) for the data $\bbf$.

\smallskip \noindent If the non-stationary subdivision scheme is convergent, and $g_{\bbf} \equiv 0$ if and only if ${\bf f}\equiv0$, then the subdivision scheme is termed \emph{non-singular}. In the forthcoming discussion we restrict ourselves to non-singular schemes only.

\smallskip \noindent
As it will be better clarified later, with respect to the subdivision capability of reproducing specific classes of functions, the standard parametrization (corresponding to the choice $t_i^{(k)}:=\frac{i}{2^k}$, $i \in \ZZ$) is not always the optimal one.
Indeed, the choice
\begin{equation}\label{def:dual_par}
   t_i^{(k)}:=\frac{i+p}{2^k},\quad i \in \ZZ, \quad p\in \RR,\quad k\ge 0\,,
\end{equation}
with $p$ suitably set, turns out to be a better selection.
In particular, when $p\in \ZZ$ the parametrization is termed \emph{primal}, whereas if $p\in \frac{\ZZ}{2}$ \emph{dual}.
For a complete discussion concerning the choice of the parametrization in the analysis of the polynomial reproduction properties of stationary subdivision schemes, we refer the reader to the papers \cite{CC13,ContiHormann11,DHSS08}.

\medskip \noindent
In consideration of the fact that the main goal of this work is the construction of a special class of non-stationary subdivision symbols capable of generating as limit functions exponential polynomials, we continue by recalling the following definitions (see, e.g, \cite{CCR13,ContiRomani11,R09}).

\begin{defn}[Exponential polynomials] \label{def:space-expol}
Let $\Gamma:=\{(\theta_1,\tau_1),\dots,(\theta_n, \tau_n)\}$ with $\theta_i\in \RR\cup \ri \RR$, $\theta_i\neq \theta_j$ if $i\neq j$ and $\tau_i\in \NN,\ i=1,\cdots,n$.
We define the space of exponential polynomials $EP_\Gamma$ as
$$
EP_\Gamma:=\hbox{span}\{\ x^{r_j} \, e^{\theta_j x},\ r_j=0,\cdots, \tau_{j}-1,\ j=1,\cdots, n\}\,.
$$
\end{defn}

\noindent For a fixed $\Gamma$, and for the corresponding space $EP_\Gamma$, we recall the following definition.

\begin{defn} [E-Generation and E-Reproduction] \label{def:ERgenerationlimit}
Let $\{a^{(k)}(z),\ k\ge 0\}$ be a sequence of subdivision symbols. The subdivision scheme associated with the symbols $\{a^{(k)}(z),\ k\ge 0\}$ is said to be \emph{$EP_\Gamma$-generating} if it is convergent and
for $f\in EP_\Gamma$ there exists an initial sequence ${\bbf}^{(0)}:=\{{\tilde f}(t^{(0)}_i),\ i\in \ZZ\}$, ${\tilde f}\in EP_\Gamma$
such that $\displaystyle{\lim_{k\rightarrow +\infty}S_{\ba^{(k)}} S_{\ba^{(k-1)}} \cdots S_{\ba^{(0)}}{\bbf}^{(0)}=f}$.
Moreover, it is said to be \emph{$EP_\Gamma$-reproducing} if it is convergent and
for $f\in EP_\Gamma$ and for the initial sequence  ${\bbf}^{(0)}:=\{f(t^{(0)}_i),\ i\in \ZZ\}$,
$\displaystyle{\lim_{k\rightarrow +\infty}S_{\ba^{(k)}} S_{\ba^{(k-1)}} \cdots S_{\ba^{(0)}}{\bbf}^{(0)}=f}$.
\end{defn}

\noindent
In the following theorem we recall the algebraic conditions on the $k$-level symbol that fully identify the generation and reproduction properties of a non-singular, univariate, binary, non-stationary subdivision scheme. A more general version of these conditions, holding for non-stationary subdivision schemes of arbitrary arity, has recently appeared in \cite{CCR13}.

\begin{theorem}\label{theo:ERreproduction}\cite[Theorem 1]{ContiRomani11}
Let $\Gamma:=\{(\theta_1,\tau_1),\dots,(\theta_n, \tau_n)\}$ with $\theta_{\ell} \in \RR\cup \ri \RR$, $\theta_{\ell}\neq \theta_j$ if $\ell \neq j$ and $\tau_{\ell}\in \NN$, $\ell=1,\cdots,n$. Let also $z^{(k)}_\ell:=e^{\frac{-\theta_\ell}{2^{k+1}}},\ \ell=1,\dots,n$. A non-singular, non-stationary subdivision scheme associated with the symbols $\{a^{(k)}(z),\ k\ge 0\}$ generates the space of exponential polynomials $EP_\Gamma$
if and only if, for each $k\ge 0$, the following conditions are satisfied
\begin{equation}\label{ERcond2}
\frac{d^r\, a^{(k)}(-z^{(k)}_\ell)}{dz^r}=0,\qquad \ell=1,\dots,n, \quad r=0,\dots,\tau_\ell-1.
\end{equation}
Furthermore, it reproduces the space of exponential polynomials $EP_\Gamma$
if and only if, for each $k\ge 0$, in addition to (\ref{ERcond2}) the following conditions are satisfied
\begin{equation}\label{ERcond1}
\frac{d^r\, a^{(k)}(z^{(k)}_\ell)}{dz^r}=2\left(z^{(k)}_\ell\right)^{p-r}\,\displaystyle{\prod_{i=0}^{r-1}(p-i)},\qquad \ell=1,\dots,n,\quad r=0,\dots,\tau_\ell-1,
\end{equation}
where an empty product is understood to be equal to $1$, and $p \in \RR$ is the shift parameter identifying the parametrization in \eqref{def:dual_par}.
\end{theorem}

\section{Symmetric subdivision symbols reproducing exponential polynomials}\label{sec:symmetry}

In this section we analyze in detail the case of $EP_{\Gamma}$-reproducing subdivision schemes featured by symmetric symbols,
since they are considered of remarkable interest in several applications.
To this purpose, we first introduce the definition of $k$-level symmetric symbol, then we point out the symmetric structure required on the set $\Gamma$ identifying the space of exponential polynomials $EP_{\Gamma}$ reproduced by a symmetric subdivision scheme.

\begin{defn}[Symmetric $k$-level symbol]\label{def_symm}
A $k$-level subdivision symbol $a^{(k)}(z)$ is called {\rm odd-symmetric} if $a^{(k)}(z)=a^{(k)}(z^{-1})$ and {\rm even-symmetric} if $z \, a^{(k)}(z)=a^{(k)}(z^{-1})$. In terms of $k$-level masks the odd/even symmetry translates into the condition $\ra^{(k)}_{-i}=\ra^{(k)}_{i}$, $i \in \ZZ$, and $\ra^{(k)}_{-i}=\ra^{(k)}_{i-1}$, $i \in \ZZ$, respectively.
\end{defn}

\begin{rem}\label{rem:shiftedsymmetry}
It is worth mentioning that a subdivision scheme has to be considered symmetric even if its $k$-level symbol satisfies the above condition after a suitable shift, \ie after multiplication by  $z^s,\ s\in \ZZ$. Note that, as shown in \cite{ContiRomani11}, the shift $s$ does affect the value of the parameter $p$ in a well-known way: the parameter $p_s$, characterizing the parametrization of the shifted scheme, is simply $p_s=p+s$.
\end{rem}

\smallskip \noindent
A symmetric set $\Gamma$ is characterized as in the following definition.

\begin{defn}[Symmetric set $\Gamma$]\label{defn:sym_gamma}
Let $\Gamma=\{\gamma_1, ..., \gamma_N\}$ be the set of cardinality $N$ in Definition \ref{def:space-expol}, containing all $\theta$ values counted with their multiplicities. The set $\Gamma$ is said to be \emph{symmetric} if
\begin{equation}\label{def:Gamma_sym}
\begin{array}{ll}
\Gamma:=\left\{
          \begin{array}{ll}
          \{(\theta_{\ell},\tau_{\ell}), (-\theta_{\ell},\tau_{\ell})\}_{\ell=1,...,\frac{N}{2}},
& \hbox{when $N$ is even,} \\ \\
           \{(\theta_{\ell},\tau_{\ell}), (-\theta_{\ell},\tau_{\ell})\}_{\ell=1,...,\frac{N-1}{2}}
\cup \{0\}, & \hbox{when $N$ is odd}
          \end{array}
        \right.
\\
\\
 \hbox{with} \quad \theta_{\ell} \in \RR^+ \cup \ri [0,\pi), \quad
\theta_{\ell} \neq \theta_j \ \ \hbox{if} \ \ \ell \neq j,
\end{array}
\end{equation}
and $\RR^+:=\{x\in \RR \, : \, x>0\}$. The space of exponential polynomials $EP_{\Gamma}$
associated to a symmetric set $\Gamma$ is also said to be \emph{symmetric}.
\end{defn}

\noindent
In the remainder of the paper we focus our attention on $EP_{\Gamma}$-reproducing symmetric subdivision schemes. We thus always assume the set $\Gamma$ to be featured by the symmetric structure specified in Definition \ref{defn:sym_gamma}.

\smallskip \noindent
The next proposition proves two very important properties of $EP_{\Gamma}$-reproducing symmetric subdivision schemes.

\begin{proposition}\label{p0p12}
A non-singular, non-stationary subdivision scheme associated with odd-symmetric or even-symmetric symbols $\{a^{(k)}(z),\ k\ge 0\}$ reproduces the pair of exponential polynomials $\{e^{\theta_{\ell} x},\ e^{-\theta_{\ell} x}\}$, $\theta_{\ell} \in \RR^+ \cup \ri (0,\pi)$, only if $p=0$  or $p=-\frac12$, respectively. Moreover, in case $\theta_{\ell}=0$, the subdivision scheme reproduces $\{1,x\}$ only if $p=0$  or $p=-\frac12$, respectively.
\end{proposition}

\proof
Let $z^{(k)}_{\ell}:=e^{\frac{-\theta_{\ell}}{2^{k+1}}}$, $\theta_{\ell} \in \RR^+ \cup \ri (0,\pi)$. We know from conditions \eqref{ERcond1} that the reproduction of
the pair $\{e^{\theta_{\ell}x}, e^{-\theta_{\ell}x} \}$ is equivalent to the existence of a shift parameter $p$ such that $a^{(k)}(z^{(k)}_{\ell})=2 (z^{(k)}_{\ell})^p$ and $a^{(k)}((z^{(k)}_{\ell})^{-1})=2 (z^{(k)}_{\ell})^{-p}$.
Thus, if  the $k$-level symbols are odd-symmetric, we can write $2 (z^{(k)}_{\ell})^p=2 (z^{(k)}_{\ell})^{-p}$, and the latter equation is satisfied only if the shift parameter $p=0$ is chosen.
Otherwise, if  the $k$-level symbols are even-symmetric, we can write $2 (z^{(k)}_{\ell})^{p+1}=2 (z^{(k)}_{\ell})^{-p}$, and the latter equation is fulfilled only if the shift parameter $p=-\frac12 $ is fixed.\\
To conclude the proof we observe that, when $\theta_{\ell}=0$, the reproduction of the pair $\{1,x\}$ is obtained by setting $p=0$ if the $k$-level symbol is odd-symmetric and $p=-\frac12$ if it is even-symmetric, as shown in \cite{ContiHormann11}.
\qed

\medskip \noindent
We continue by analyzing useful algebraic properties fulfilled by symmetric subdivision symbols.

\begin{proposition}\label{evensym_yoon}
Let $\Gamma=\{(\theta_{\ell}, d), (-\theta_{\ell}, d)\}$ with $\theta_{\ell} \in \RR^+ \cup \ri [0,\pi)$ be a set of cardinality $2d$, $d \in \NN$.
For $z^{(k)}_{\ell}=e^{\frac{-\theta_{\ell}}{2^{k+1}}}$ the even-symmetric subdivision symbols $\{a^{(k)}(z),\ k\ge 0\}$ satisfy
$$\frac{d^r\, a^{(k)}(z^{(k)}_{\ell})}{dz^r}=2\left(z^{(k)}_{\ell} \right)^{-\frac12-r}\,\displaystyle{\prod_{i=0}^{r-1}\left(-\frac{1}{2}-i\right)},\quad r=0,\dots,d-1,$$
if and only if  the associated odd-symmetric subdivision symbols $\{b^{(k)}(z), \ k \ge 0\}$ with $b^{(k)}(z)=z \, a^{(k)}(z^2) -2$ satisfy
$$
\frac{d^r b^{(k)}((z^{(k)}_{\ell})^{\frac12})}{dz^r}=0, \qquad  r=0,\dots,d-1. $$ In the above equation, an empty product is understood to be equal to $1$.
\end{proposition}

\proof
We start showing that the $k$-level symbol
$a^{(k)}(z)$ is even-symmetric if and only if
$b^{(k)}(z)=z \, a^{(k)}(z^2) -2$ is odd-symmetric. Indeed $z \, a^{(k)}(z)=a^{(k)}(z^{-1})$ if and only if
$z^2 \, a^{(k)}(z^2)=a^{(k)}(z^{-2})$ if and only if
$b^{(k)}(z)+2=b^{(k)}(z^{-1})+2$ if and only if $b^{(k)}(z)=b^{(k)}(z^{-1})$.\\
The rest of the proof is inductive on $r$. The case $r=0$ is easy to check.
Therefore we consider
the case $r>0$ and use the Leibniz formula and the induction for $r=0$ to write the derivatives of $a^{(k)}(z)=z^{-\frac12}(b^{(k)}(z^{\frac12})+2)$ evaluated at $z^{(k)}_{\ell}$. Recall that
$z^{-\frac12}=e^{-\frac12 {\rm log}(z)}$  can be defined as  a single-valued  function,
analytic on $\mathbb C\setminus (-\infty, 0]$. Thus we have
\[
\begin{array}{lll}
\frac{d^r a^{(k)}(z)}{dz^r} \Big \vert_{z=z^{(k)}_{\ell}}=\sum_{s=0}^r \binom{r}{s} \frac{d^s (b^{(k)}(z^{\frac12})+2)}{dz^s} \frac{d^{r-s} (z^{-\frac12})}{dz^{r-s}} \Big \vert_{z=z^{(k)}_{\ell}}\\
                                                                     \\
                                                                     =\left(\frac{d^{r} (z^{-\frac12})}{dz^{r}}(b^{(k)}(z^{\frac12})+2)\right)
                                                                     \Big \vert_{z=z^{(k)}_{\ell}}+ \sum_{s=1}^r
                                                                     \binom{r}{s}
                                                                     \frac{d^s (b^{(k)}(z^{\frac12})+2)}{dz^s}
                                                                     \frac{d^{r-s} (z^{-\frac12})}{dz^{r-s}} \Big \vert_{z=z^{(k)}_{\ell}}\\
                                                                     \\
                                                                     =2\prod_{i=0}^{r-1} \left(-\frac12-i \right)(z^{(k)}_{\ell})^{-\frac12-r}
+  \sum_{s=1}^r \binom{r}{s}
  \frac{d^s (b^{(k)}(z^{\frac12}))}{dz^s} \frac{d^{r-s} (z^{-\frac12})}{dz^{r-s}} \Big \vert_{z=z^{(k)}_{\ell}}.
 \end{array}
\]
 We continue by using the Fa\`a di Bruno's formula (see \cite{WJ02} or \cite{Mortini}) to write
 $$
 \frac{d^s (b^{(k)}(z^{\frac12}))}{dz^s} \Big \vert_{z=z^{(k)}_{\ell}} =
\sum_{j=1}^s \frac{d^j b^{(k)}(y)}{dy^j} \Big \vert_{y=(z^{(k)}_{\ell})^\frac12} \, A_{j,s}(z) \Big \vert_{z=z^{(k)}_{\ell}},
$$
with $A_{j,s}(z)$
%=\sum_{\bq \in \bM^j, \, \vert \bq \vert=r} \frac{C^r_j(z)}{\prod_{i=1}^r N(\bq,i)!}\,$
given functions whose value is important to know only for $j=s$. In particular, we have $A_{s,s}(z)=\left(\frac12 z^{-\frac12}\right)^s$. In fact, for $j=1,\cdots,s,\ s<r$, using the induction assumption  we know that $\frac{d^j b^{(k)}(y)}{dy^j} \Big \vert_{y=(z^{(k)}_{\ell})^\frac12}=0$ and therefore the above sum reduces to the last term only, that is to
 $$ \frac{d^s b^{(k)}(y)}{dy^s} \Big \vert_{y=(z^{(k)}_{\ell})^\frac12} \left(\frac12 (z^{(k)}_{\ell})^{-\frac12}\right)^s.$$
 Hence, using the fact that $\frac{d^r\, a^{(k)}(z)}{dz^r} \Big \vert_{z=z_{\ell}^{(k)}}=2\left(z^{(k)}_{\ell}\right)^{-\frac12-r}\prod_{i=0}^{r-1} \left(-\frac12-i \right)$, we arrive at
\[
 \begin{array}{ll}
2\left(z^{(k)}_{\ell}\right)^{-\frac12-r}\,\prod_{i=0}^{r-1} \left(-\frac12-i \right)=2{\prod_{i=0}^{r-1} \left(-\frac12-i \right)\left(-\frac12\right)}(z^{(k)}_{\ell})^{-\frac12-r}\\
+
\sum_{s=1}^r \binom{r}{s}
\frac{d^s b^{(k)}(y)}{dy^s} \Big \vert_{y=(z^{(k)}_{\ell})^\frac12} \left(\frac12 (z^{(k)}_{\ell})^{-\frac12}\right)^s \,\prod_{i=0}^{r-s-1} \left(-\frac12-i \right) \left(-\frac12 \right) (z_{\ell}^{(k)})^{-\frac12-(r-s)}.
 \end{array}
\]
 Now, using again the inductive hypothesis that $\frac{d^s b^{(k)}(y)}{dy^s} \Big \vert_{y=(z^{(k)}_{\ell})^\frac12}=0$ for all $s=1,...,r-1$ we obtain
 $$
0=
\frac{d^r b^{(k)}(y)}{dy^r} \Big \vert_{y=(z^{(k)}_{\ell})^\frac12} \, \left(\frac12 (z^{(k)}_{\ell})^{-\frac12}\right)^r \, (z_{\ell}^{(k)})^{-\frac12},
 $$
which is the required value of the $r$-th derivative of $b^{(k)}(z)$ at $(z^{(k)}_{\ell})^\frac12$, i.e.
 $$
 \frac{d^r b^{(k)}(y)}{dy^r} \Big \vert_{y=(z^{(k)}_{\ell})^\frac12} =0\,.
 $$
 This concludes the induction step and therefore the proof.
\qed

\begin{rem}\label{add1}
As a by-product of the proof of Proposition \ref{evensym_yoon}  from $b^{(k)}(z)=b^{(k)}(z^{-1})$
we obtain that, for $z^{(k)}_{\ell}=1$, the condition
$\frac{d^r b^{(k)}((z^{(k)}_{\ell})^{\frac12})}{dz^r}=0$  for $r=0,\dots,2j$,  with
$0\leq 2j\leq \tau_{\ell}-1$,   implies that
 $\frac{d^{2j+1} b^{(k)}((z^{(k)}_{\ell})^{\frac12})}{dz^{2j+1}}=0$ and, therefore, the corresponding condition on
$\frac{d^{2j+1} a^{(k)}(z^{(k)}_{\ell})}{dz^{2j+1}}$.
\end{rem}

\noindent
The next proposition shows that conditions \eqref{ERcond1} are compatible with symmetry properties of subdivision symbols. Indeed we prove that symmetric subdivision symbols are such that, if conditions \eqref{ERcond1} are satisfied at a given $z_{\ell}^{(k)}$, they are also satisfied at $(z_{\ell}^{(k)})^{-1}$.

\begin{proposition}\label{prop:oddsymm_EPR}
Let $\{a^{(k)}(z),\ k\ge 0\}$ be the odd-symmetric (even-symmetric) symbols of a non-singular, non-stationary subdivision scheme associated with the shift parameter $p=0$ ($p=-\frac12$). For $z^{(k)}_{\ell}:=e^{\frac{-\theta_{\ell}}{2^{k+1}}}$ with $\theta_{\ell} \in \RR^+ \cup \ri [0,\pi)$  we have that
$$\frac{d^r\, a^{(k)}(z^{(k)}_{\ell})}{dz^r}=2\left(z^{(k)}_{\ell}\right)^{p-r}\,\prod_{i=0}^{r-1}(p-i),\quad r=0,\ldots, d-1,\quad  d\in \NN$$
if and only if
$$\frac{d^r\, a^{(k)}((z^{(k)}_{\ell})^{-1})}{dz^r}=2\left((z^{(k)}_{\ell})^{-1}\right)^{p-r}\,\prod_{i=0}^{r-1}(p-i)\, \quad r=0,\ldots, d-1,\quad  d\in \NN\,,$$
where an empty product is understood to be equal to $1$.
\end{proposition}

\proof
We show the claim by induction on $r$. The case $r=0$ has been already considered in Proposition \ref{p0p12}.
 For $r>0$ we first consider the odd-symmetric case and we start proving one of the two implications. Computing the $r$-th derivative of the equation $a^{(k)}(z)=a^{(k)}(z^{-1})$ via the Fa\`a di Bruno's formula (see \cite{WJ02} or \cite{Mortini}) and evaluating it at $z^{(k)}_{\ell}$, we obtain
$$
\frac{d^r a^{(k)}(z)}{dz^r} \Big \vert_{z=z^{(k)}_{\ell}}= \sum_{j=1}^r \frac{d^j a^{(k)}(y)}{dy^j} \Big \vert_{y=(z^{(k)}_{\ell})^{-1}} A_{j,r}(z^{(k)}_{\ell})
$$
with $$A_{j,r}(z)=\sum_{\bq \in \bM^j, \, \vert \bq \vert=r} \frac{r!}{\bq!} \frac{(-1)^r z^{-r-j}}{\prod_{i=1}^r N(\bq,i)!}\,,
$$
where
$$
\bM^j=\{\bq=(q_1, q_2, ...,q_j) \in \NN^j, \, q_1 \geq q_2 \geq ... \geq q_j \geq 1  \}, \qquad \vert \bq \vert=q_1+...+q_j,
$$
and $N(\bq,i)$ denoting the number of times the positive integer $i$ appears in the $j$-tuple $\bq \in \NN^j$.
Now, by the inductive hypothesis we know that $\frac{d^r a^{(k)}(z)}{dz^r} \Big \vert_{z=z^{(k)}_{\ell}}=0$ for $r=1,...,d-2$ implies $\frac{d^r a^{(k)}(z)}{dz^r} \Big \vert_{z=(z^{(k)}_{\ell})^{-1}}=0$ for $r=1,...,d-2$. Hence,
$$
\frac{d^{d-1} a^{(k)}(z)}{dz^{d-1}} \Big \vert_{z=z^{(k)}_{\ell}}=  \frac{d^{d-1} a^{(k)}(y)}{dy^{d-1}} \Big \vert_{y=(z^{(k)}_{\ell})^{-1}} A_{d-1,d-1}(z^{(k)}_{\ell}),
$$
and since $A_{d-1,d-1}(z^{(k)}_{\ell})=(-1)^{d-1}(z^{(k)}_{\ell})^{-2(d-1)}\neq 0$, we easily get that
$$\frac{d^{d-1} a^{(k)}(z)}{dz^{d-1}} \Big \vert_{z=z^{(k)}_{\ell}}=0 \quad \Rightarrow \quad \frac{d^{d-1} a^{(k)}(y)}{dy^{d-1}} \Big \vert_{y=(z^{(k)}_{\ell})^{-1}}=0\,,$$
which concludes one direction of the proof in the odd-symmetric case. The proof of the converse implication can be repeated analogously.

\smallskip \noindent For the even-symmetric case, in view of Proposition \ref{evensym_yoon}, we use the same argument as above for the odd-symmetric $k$-level symbol
$b^{(k)}(z)=z \, a^{(k)}(z^2) -2$ and for the roots $(z^{(k)}_{\ell})^\frac12,\ (z^{(k)}_{\ell})^{-\frac12}$, so completing the proof.
\qed

\smallskip
\begin{rem}
The above proposition proves the equivalence between the conditions for exponential polynomial reproduction given in \cite{ContiRomani11} and \cite{YoonAD2013,YoonMA2013} when $p=0$ or $p=-\frac12$.
\end{rem}

\section{Symmetric exponential B-splines and their normalization factors}\label{exponential B-splines}

For a symmetric set  $\Gamma$ as in Definition \ref{defn:sym_gamma}, we introduce the following notation where, for a given $L \in \RR$, $\lfloor L  \rfloor:=\max \{M \in \ZZ \, : M \leq L\}$.
For $\ell=1,\cdots, \lfloor \frac N2 \rfloor$, in case $N$ is even we denote by $\Gamma_{\ell,e}:=\Gamma \setminus\{\theta_\ell,-\theta_\ell\}$, and in case $N$ is odd we define $\Gamma_{\ell,o}:=\Gamma \setminus \{\theta_\ell,-\theta_\ell,0\}$.\\
In the following, for a given $L \in \RR$, we also use the notation $\lceil L  \rceil:=\min \{M \in \ZZ \, : M \geq L\}$.

\medskip \noindent
A symmetric (not-normalized) exponential B-spline is based on the sequence of symbols
\begin{equation}\label{def:nnES}
{\Bt}_{N,\Gamma}^{(k)}(z):=z^{-\lceil \frac N2 \rceil}\prod_{i=1}^N\left(
e^{\frac{\gamma_i}{2^{k+1}}}z+1\right),\quad
\gamma_i \in \Gamma, \qquad
k\ge 0\,.
\end{equation}
By definition of $\Gamma_{\ell,e}$ we easily see that $\Bt_{N,\Gamma}^{(k)}(z)$ satisfies a ``recursion'' formula since

\begin{equation}\label{def:recursion}
{\Bt}_{N,\Gamma}^{(k)}(z)=z^{-1}\left(
e^{\frac{\theta_\ell}{2^{k+1}}}z+1\right)\left(
e^{\frac{-\theta_\ell}{2^{k+1}}}z+1\right){\Bt}_{N-2,\Gamma_{\ell,e}}^{(k)}(z)\,, \qquad \ell=1,\dots, \left \lfloor\frac{N}{2} \right\rfloor.
\end{equation}
For later use we also observe that for any $\theta_i \neq |\theta_\ell|$
\begin{equation}\label{prop:Bevaluation}
\left\{
  \begin{array}{lll}
  &  {\Bt}_{N-2,\Gamma_{\ell,e}}^{(k)} \left(e^{\frac{\theta_i}{2^{k+1}}}\right)={\Bt}_{N-2,\Gamma_{\ell,e}}^{(k)} \left(e^{\frac{-\theta_i}{2^{k+1}}}\right)\,, & \hbox{if $N$ is even} \smallskip\\
  &  (e^{\frac{-\theta_i}{2^{k+1}}}+1){\Bt}_{N-2,\Gamma_{\ell,o}}^{(k)} \left(e^{\frac{\theta_i}{2^{k+1}}} \right)= (e^{\frac{\theta_i}{2^{k+1}}}+1){\Bt}_{N-2,\Gamma_{\ell,o}}^{(k)} \left(e^{\frac{-\theta_i}{2^{k+1}}}\right)\,, & \hbox{if $N$ is odd.}
  \end{array}
\right.
\end{equation}
Moreover, the symbols in (\ref{def:nnES}) satisfy the necessary and sufficient conditions for $EP_{\Gamma}$-generation
\begin{equation}\label{prop:Egeneration}
\left\{
  \begin{array}{ll}
\begin{array}{lll}
   {\Bt}_{N,\Gamma}^{(k)}\left(-e^{\frac{\pm \theta_j}{2^{k+1}}}\right)=0,\\
\displaystyle \frac{d^{r}\, {\Bt}_{N,\Gamma}^{(k)}\left(-e^{\frac{\pm \theta_j}{2^{k+1}}}\right)}{dz^{r}}=0,
\quad r=1,\dots,\tau_j-1,\ & j=1,\dots,\frac{N}{2},\\
  \hbox{if $N$ is even};
\end{array}
\\
\begin{array}{lll}
   {\Bt}_{N,\Gamma}^{(k)}(-1)=0, \quad  {\Bt}_{N,\Gamma}^{(k)}\left(-e^{\frac{\pm \theta_j}{2^{k+1}}}\right)=0,\\
 \displaystyle\frac{d^{r}\, {\Bt}_{N,\Gamma}^{(k)}\left(-e^{\frac{\pm \theta_j}{2^{k+1}}}\right)}{dz^{r}}=0,
\quad r=1,\dots,\tau_j-1,\ & j=1,\dots,\frac{N-1}{2},\\
 \hbox{if $N$ is odd.}
\end{array}
  \end{array}
\right.
\end{equation}

\smallskip \noindent For reproduction purposes it may be convenient to consider \emph{normalized exponential B-splines}. Their symbols are defined by multiplying ${\Bt}_{N,\Gamma}^{(k)}(z)$ in (\ref{def:nnES}) with an extra factor $K_{\ell}^{(k)} \in \RR$, namely by

\begin{equation}\label{def:nES}
B_{N,\Gamma}^{(k)}(z):=K_{\ell}^{(k)}\, {\Bt}_{N,\Gamma}^{(k)}(z)=K_{\ell}^{(k)} \, z^{-\lceil \frac N2 \rceil}\prod_{i=1}^N\left(
e^{\frac{\gamma_i}{2^{k+1}}}z+1\right),\quad
\gamma_i \in \Gamma, \qquad
k\ge 0\,,
\end{equation}
where the $k$-level coefficient $K_{\ell}^{(k)}$ can be selected in accordance with the  parameter $p$ in order to ensure that the \emph{normalized exponential B-spline}, besides generating $EP_{\Gamma}$,  \emph{reproduces} the pair of exponential polynomials $\{e^{\theta_{\ell} \, x},\ e^{-\theta_{\ell }\, x}\} \in EP_{\Gamma}$. This fact is discussed in the next proposition.

\begin{proposition}\label{prop:EreprodK}
Let $\Gamma$ be given and $\theta_{\ell} \in \Gamma$. The symbols in (\ref{def:nES}) satisfy the $\{e^{\theta_{\ell} \, x},\ e^{-\theta_{\ell }\, x}\}$-reproduction condition if
\begin{itemize}
\item[(i)] $p=0$\ \hbox{and}\  $(K_{\ell}^{(k)})^{-1}=\left(e^{\frac{-\theta_{\ell}}{2^{k+1}}}+e^{\frac{\theta_{\ell}}{2^{k+1}}}\right) \, {\Bt}_{N-2,{ \Gamma_{\ell,e}}}^{(k)}(e^{\frac{\theta_{\ell}}{2^{k+1}}})$, \ for $N$ even,
\item[(ii)] $p=-\frac12$\ \hbox{and}\  $(K_{\ell}^{(k)})^{-1}=\left(e^{\frac{-\theta_{\ell}}{2^{k+2}}}+e^{\frac{\theta_{\ell}}{2^{k+2}}}\right) \, \left(e^{\frac{-\theta_{\ell}}{2^{k+1}}}+e^{\frac{\theta_{\ell}}{2^{k+1}}}\right) \, {\Bt}_{N-3,{\Gamma_{\ell,o}}}^{(k)}(e^{\frac{\theta_{\ell}}{2^{k+1}}})$,  \ for $N$ odd.
  \end{itemize}
\end{proposition}

\proof
Let us start analyzing the case $N$ even. Introducing the abbreviation $z_{\ell}^{(k)}:=e^{\frac{-\theta_{\ell}}{2^{k+1}}}$, in view of Theorem \ref{theo:ERreproduction} the reproduction of $\{e^{\theta_{\ell} \, x},\ e^{-\theta_{\ell }\, x}\}$ requires the fulfillment of the conditions
$$
B^{(k)}_{N,\Gamma}( z_\ell^{(k)})=2 (z_\ell^{(k)})^p, \qquad  B^{(k)}_{N,\Gamma}\left((z_\ell^{(k)})^{-1}\right)=2 (z_\ell^{(k)})^{-p},
$$
that is
$$
K_{\ell}^{(k)} \, \Bt^{(k)}_{N,\Gamma}\left(z_\ell^{(k)}\right)=2 (z_\ell^{(k)})^{p}, \qquad
K_{\ell}^{(k)} \, \Bt^{(k)}_{N,\Gamma}\left((z_\ell^{(k)})^{-1}\right)=2 \left(z_\ell^{(k)}\right)^{-p}.
$$
Exploiting the recurrence relation in (\ref{def:recursion}) and recalling that
$ \Bt^{(k)}_{N-2,\Gamma_{\ell,e}}\left( z_\ell^{(k)}\right) = \Bt^{(k)}_{N-2,\Gamma_{\ell,e}}\left( (z_\ell^{(k)})^{-1}\right)$, we have
$$
\begin{array}{ll}
K_{\ell}^{(k)} \, \left( ( z_\ell^{(k)})^{-2}+1\right)  \, \Bt^{(k)}_{N-2,\Gamma_{\ell,e}}\left((z_\ell^{(k)})^{-1}\right) = (z_\ell^{(k)})^{p-1},
 \\
K_{\ell}^{(k)} \, \left( ( z_\ell^{(k)})^{-2} +1\right) \, \Bt^{(k)}_{N-2,\Gamma_{\ell,e}}\left((z_\ell^{(k)})^{-1}\right) = (z_\ell^{(k)})^{-p-1}.
\end{array}
$$
The solution of this system in the unknowns $p$ and $K_{\ell}^{(k)}$ is therefore given by
$$p=0\quad \hbox{and}\quad (K_{\ell}^{(k)})^{-1} = \left( z_\ell^{(k)}+ (z_\ell^{(k)})^{-1} \right) \Bt^{(k)}_{N-2,\Gamma_{\ell,e}}\left((z_\ell^{(k)})^{-1}\right), $$
which concludes the proof of subcase $(i)$.\\
We continue studying the case $N$ odd. Again, using the recurrence relation (\ref{def:recursion}), the two conditions to be satisfied for the reproduction of $\{e^{\theta_{\ell} \, x},\ e^{-\theta_{\ell }\, x}\}$ can be written as
$$
2 K_{\ell}^{(k)} \left(z_\ell^{(k)}\right)^{-2} \left( (z_\ell^{(k)})^2+1 \right) \left( z_\ell^{(k)}+1 \right) \Bt^{(k)}_{N-3,\Gamma_{\ell,o}}
\left( z_\ell^{(k)}\right)=2 (z_\ell^{(k)})^{p}, $$
$$2 K_{\ell}^{(k)} \left(z_\ell^{(k)}\right)^{2} \left( (z_\ell^{(k)})^{-2}+1\right) \left((z_\ell^{(k)})^{-1}+1 \right) \Bt^{(k)}_{N-3,\Gamma_{\ell,o}}\left((z_\ell^{(k)})^{-1}\right)=2 (z_\ell^{(k)})^{-p}.
$$
Now, since $N-3$ is even and
$\Bt^{(k)}_{N-3,\Gamma_{\ell,o}}\left((z_\ell^{(k)})^{-1}\right)=\Bt^{(k)}_{N-3,\Gamma_{\ell,o}
}\left( z_\ell^{(k)}\right)$, we can write the simplified expressions
$$
K_{\ell}^{(k)} \, \left( z_\ell^{(k)}+1 \right) \left((z_\ell^{(k)})^{-1}+z_\ell^{(k)}\right)  \, \Bt^{(k)}_{N-3,\Gamma_{\ell,o}}\left( (z_\ell^{(k)})^{-1}\right)= (z_\ell^{(k)})^{p+1},
$$
$$
K_{\ell}^{(k)} \, \left( z_\ell^{(k)}+1 \right) \left( (z_\ell^{(k)})^{-1} +z_\ell^{(k)}\right) \, \Bt^{(k)}_{N-3,\Gamma_{\ell,o}}\left( (z_\ell^{(k)})^{-1}\right)= (z_\ell^{(k)})^{-p}.
$$
The solution of this system in the unknowns $p$ and $K_{\ell}^{(k)}$ is given by
$$p=-\frac{1}{2} \; \hbox{and} \;
(K_{\ell}^{(k)})^{-1} = \left( (z_\ell^{(k)})^{\frac{1}{2}}+ (z_\ell^{(k)})^{-\frac{1}{2}} \right) \left( z_\ell^{(k)}+(z_\ell^{(k)})^{-1} \right) \, \Bt^{(k)}_{N-3,\Gamma_{\ell,o}}\left((z_\ell^{(k)})^{-1}\right),
$$
which concludes the proof of subcase $(ii)$.
\qed

\medskip \noindent
Note that similar results concerning the normalization of exponential B-splines  are also given in \cite{YoonAD2013}.
Two special situations are considered in the next result.

\begin{coro}
For $N\ge 2$ the exponential B-splines $B_{N, \Gamma} ^{(k)}(z)$ reproduce $\{1,\ x\}$ if $\ 0=\theta_{\ell}\in \Gamma$ with $\tau_\ell=2$,
$$
(K_{\ell}^{(k)})^{-1}=\left\{
                                              \begin{array}{ll}
                                                2 \Bt^{(k)}_{N-2,\Gamma  \setminus \{\theta_\ell, -\theta_\ell\}}(1), & \hbox{for $N$ even;} \\
                                                4 \Bt^{(k)}_{N-3,\Gamma  \setminus \{\theta_\ell, -\theta_\ell, 0\}}(1),  & \hbox{for $N$ odd;}
                                              \end{array}
                                            \right.
$$
and
$$ p=\left\{
                                              \begin{array}{ll}
                                                0, & \hbox{for $N$ even;} \\
                                                -\frac12, & \hbox{for $N$ odd.}
                                              \end{array}
                                            \right.
$$
Moreover, when $\Gamma=\{(\theta_1,\tau_1)=(0,N)\}$, the symbol of $B_{N, \Gamma } ^{(k)}(z)$ does not depend on $k$ any longer and becomes the stationary symbol of the shifted order-$N$ (polynomial) B-spline that we simply denote by
$
B_{N}(z)=z^{-\lceil\frac{N}{2}\rceil} \, \frac{(1+z)^{N}}{2^{N-1}}\,.
$

\end{coro}

\section{Deriving the symbols of exponential pseudo-splines and investigating their symmetry properties}\label{L-exponential pseudo-splines}
For any $p\in\RR$ and $N,M\in\NN$ the binary pseudo-spline subdivision scheme is defined to be the stationary scheme with minimal support that generates polynomials of degree $N-1$ and whose symbol, $a_{M,N}(z)$, satisfies the conditions $\frac{d^r a_{M,N}(1)}{dz^r}  = 2 \prod_{i=0}^{r-1}(p-i)$, $r=0,\dots,M-1$ for reproduction of polynomials up to degree $M-1$. Its actual degree of polynomial reproduction is thus $\min\{N-1, M-1\}$ (see \cite{ContiHormann11,DongShen2006a,DongShen2006b}).

\smallskip \noindent
The main contribution of this paper consists in generalizing the family of binary pseudo-splines to the non-stationary setting. The resulting family is called the family of \emph{binary exponential pseudo-splines}.

\medskip \noindent For $\Gamma$ as in \eqref{def:Gamma_sym} and $p=0$ if $N$ is even, $p=-\frac12$ if $N$ is odd,
the family of symmetric binary exponential pseudo-splines is defined to be the family of symmetric subdivision schemes with \emph{minimal support} that generates the space of exponential polynomials $EP_{\Gamma}$ and whose $k$-level symbol satisfies the conditions in \eqref{ERcond1}
for reproduction of elements in $EP_{\tilde \Gamma}$, where
${\tilde \Gamma}$ denotes a symmetric subset of $\Gamma$ of cardinality $M \leq N$, with $M$ and $N$ of the same parity.

\smallskip \noindent
The $k$-level symbol of an exponential pseudo-spline is therefore of the form
\[
  a_{M,N,\Gamma}^{(k)}(z) := B^{(k)}_{N,\Gamma}(z)\, c_{M,\Gamma}^{(k)}(z),
\]
where $B^{(k)}_{N,\Gamma}(z)$ is the $k$-level symbol of the normalized exponential B-spline in (\ref{def:nES})
with $K_{\ell}^{(k)}$ as in Proposition \ref{prop:EreprodK},
whereas $c_{M,\Gamma}^{(k)}(z)$ is the $k$-level Laurent polynomial of lowest possible degree such that $a_{M,N,\Gamma}^{(k)}(z)$ satisfies
\begin{equation}\label{eq:conditions}
\left\{\begin{array}{lll}
      a_{M,N,\Gamma}^{(k)}(z^{(k)}_\ell)=2 \, \left(z^{(k)}_{\ell}\right)^p,
\\
 \displaystyle\frac{d^s\, a_{M,N,\Gamma}^{(k)}(z^{(k)}_\ell)}{dz^s}=2\left(z^{(k)}_\ell\right)^{p-s}\,\displaystyle{\prod_{i=0}^{s-1}(p-i)},
\\
\ell=1,\dots,m, \quad m \leq n; \quad s=1,\dots,\tau_\ell-1,
 \end{array}\right.
\end{equation}
for $z^{(k)}_\ell:=e^{\frac{-\theta_\ell}{2^{k+1}}}$ and $(\theta_\ell,\tau_\ell)\in {\tilde \Gamma} \subset \Gamma$ \, for $\ell=1,...,m$, with $M:=\sum_{j=1}^m \tau_j$.
Obviously, (\ref{ERcond2}) are satisfied by construction.

\medskip \noindent Using the Leibniz rule we can write the set of conditions in (\ref{eq:conditions}) in the equivalent form
\begin{equation}\label{eq:Leibniz}
  \sum_{i=0}^s \binom{s}{i}\frac{d^i\, c_{M,\Gamma}^{(k)}(z^{(k)}_\ell)}{dz^{i}} \frac{d^{s-i}\, B^{(k)}_{N,\Gamma}(z^{(k)}_\ell)}{dz^{s-i}}=
  v_{\ell,s},\quad
      \begin{array}{l}
      \ell=1,\dots,m, \quad m \leq n\\
     s=0,\dots,\tau_\ell-1,
\end{array}
\end{equation}
where $v_{\ell,s}:=2\left(z^{(k)}_\ell\right)^{p-s}\,\displaystyle{\prod_{i=0}^{s-1}(p-i)}$,
and an empty product is understood to be equal to $1$.

\medskip \noindent
We start by considering the case $m=n$ (which means $M=N$). In the latter case equations \eqref{eq:Leibniz} can be rewritten as the linear system
\begin{equation}\label{eq:linear-system-pseudo-splines}
  \bcalA \bw= \bv,\quad  \bcalA \in \RR^{N\times N},\quad \bw\in \RR^{N},\quad  \bv\in \RR^{N},\quad N=\sum_{j=1}^n\tau_j,
\end{equation}
where $\bw:=\Bigl( \frac{d^i\, c_{N,\Gamma}^{(k)}(z^{(k)}_\ell)}{dz^{i}},\ i=0,\cdots,\tau_\ell-1,\ \ell=1,\cdots,n\Bigr)^T$,
$\bv$ is defined as $\bv:=(v_{1,0},\ldots, v_{1, \tau_1-1},v_{2,0}, \ldots, v_{2, \tau_2-1},  \ldots)^T$,
and $\bcalA$ is the block diagonal lower triangular matrix  given by
\[
\bcalA :=\left[\begin{array}{lll}
\mathcal A_1 \\ &  \ddots \\ && \mathcal A_n \end{array}\right], \quad \mathcal A_j\in \mathbb C^{\tau_j\times \tau_j},
\]
where,  for $1\leq j\leq n$,
\[
\mathcal A_j :=\left[\begin{array}{cccc}
B^{(k)}_{N,\Gamma}(z^{(k)}_j) & & & \smallskip \\
\binom{1}{0}\frac{d^{1}\, B^{(k)}_{N,\Gamma}(z^{(k)}_j)}{dz^{1}} & B^{(k)}_{N,\Gamma}(z^{(k)}_j) & & \smallskip\\
\vdots & & \ddots & \smallskip \\
\binom{\tau_j-1}{0}\frac{d^{\tau_{j}-1}\, B^{(k)}_{N,\Gamma}(z^{(k)}_j)}{dz^{\tau_j-1}}&
\binom{\tau_j-1}{1}\frac{d^{\tau_{j}-2}\, B^{(k)}_{N,\Gamma}(z^{(k)}_j)}{dz^{\tau_j-2}}& \ldots &
B^{(k)}_{N,\Gamma}(z^{(k)}_j)\end{array}\right].
\]

\smallskip \noindent
The structure of $\mathcal A_j^{-1}$ follows from the next  result.
\begin{lemma}\label{lem:inverse}
For a given $z^{(k)}_j$ such that $B^{(k)}_{N,\Gamma}(z^{(k)}_j)\neq 0$ the matrix $\mathcal A_j$ is invertible and
\[
\mathcal A_j^{-1}=\mathcal G_j :=
\left[\begin{array}{cccc}
G^{(k)}_{N,\Gamma}(z^{(k)}_j) \smallskip\\
\binom{1}{0}\frac{d^{1}\, G^{(k)}_{N,\Gamma}(z^{(k)}_j)}{dz^{1}} & G^{(k)}_{N,\Gamma}(z^{(k)}_j) \smallskip\\
\vdots & & \ddots \smallskip\\
\binom{\tau_j-1}{0}\frac{d^{\tau_j-1}\, G^{(k)}_{N,\Gamma}(z^{(k)}_j)}{dz^{\tau_j-1}}&
\binom{\tau_j-1}{1}\frac{d^{\tau_j-2}\, G^{(k)}_{N,\Gamma}(z^{(k)}_j)}{dz^{\tau_j-2}}& \ldots &
G^{(k)}_{N,\Gamma}(z^{(k)}_j)\end{array}\right]
\]
with $G^{(k)}_{N,\Gamma}(z):=1/B^{(k)}_{N,\Gamma}(z)$.
\end{lemma}

\proof
Let $\mathcal S_j=\mathcal A_j \, \mathcal G_j$ be the product matrix.
The claim follows from the relation $B^{(k)}_{N,\Gamma}(z) \ G^{(k)}_{N,\Gamma}(z)=1$.
By differentiating and using the Leibniz rule we obtain that for $s\geq 1$
\[
\sum_{i=0}^s \binom{s}{i}\frac{d^i\,G^{(k)}_{N,\Gamma}(z^{(k)}_j)}{dz^{i}} \frac{d^{s-i}\, B^{(k)}_{N,\Gamma}(z^{(k)}_j)}{dz^{s-i}}=0,
\]
which means that the subdiagonal entries in the first column of $\mathcal S_j$ are zero.  For the remaining
entries observe that
\[
(\mathcal S_j)_{\ell,r}=
\sum_{i=r}^\ell \binom{\ell-1}{i-1}\binom{i-1}{r-1}\frac{d^{\ell-i}\, B^{(k)}_{N,\Gamma}(z^{(k)}_j)}{dz^{\ell-i}}
\frac{d^{i-r}\, G^{(k)}_{N,\Gamma}(z^{(k)}_j)}{dz^{i-r}},
\]
and, hence, by setting $\tilde i := i-r$ and $\tilde \ell := \ell-r$
\[
(\mathcal S_j)_{\ell,r}=\binom{\ell-1}{r-1}
\sum_{\tilde i=0}^{\tilde \ell} \binom{\tilde \ell}{\tilde i}\frac{d^{\tilde \ell-\tilde i}\, B^{(k)}_{N,\Gamma}(z^{(k)}_j)}
{dz^{\tilde \ell-\tilde i}}
\frac{d^{\tilde i}\, G^{(k)}_{N,\Gamma}(z^{(k)}_j)}{dz^{\tilde i}}=\delta_{\ell,r},
\]
where $\delta_{\ell,r}$ denotes the Kronecker symbol.
\qed

\smallskip \noindent
Assuming that $\Gamma$ is defined as in \eqref{def:Gamma_sym}, we have that  $z^{(k)}_j\neq -z^{(k)}_\ell$, $1\leq j, \ell\leq N$ and thus Lemma \ref{lem:inverse}  yields  the following.

\begin{proposition}\label{prop:derivateC}
For $z^{(k)}_\ell:=e^{\frac{-\theta_\ell}{2^{k+1}}},\ \ell=1,\cdots,n$, we have
\begin{equation}\label{intgen}
\frac{d^s\, c_{N,\Gamma}^{(k)}(z^{(k)}_\ell)}{dz^{s}}=
\sum_{i=0}^s \binom{s}{i}  v_{\ell, i}\frac{d^{s-i}\, G^{(k)}_{N,\Gamma}(z^{(k)}_\ell)}{dz^{s-i}}, \quad
 \begin{array}{l}
      \ell=1,\dots,n, \\
      s=0,\dots,\tau_\ell-1,
\end{array}
\end{equation}
where $v_{\ell, i}:=2\left(z^{(k)}_\ell\right)^{p-i}\,\displaystyle{\prod_{r=0}^{i-1}(p-r)}$ and  $G^{(k)}_{N,\Gamma}(z):=\frac 1{B^{(k)}_{N,\Gamma}(z)}$.
\end{proposition}

\smallskip\noindent
It is worth pointing out that, once we have specified the support of $c_{N,\Gamma}^{(k)}(z)$, the generalized interpolation conditions
\eqref{intgen}  enable the computation of  its coefficients by means of an interpolation process.
In particular, the construction of $c_{N,\Gamma}^{(k)}(z)$ can rely upon the following functional approach.  Let
$\{\beta_i\}_{i=1}^N$, $N=\sum_i\tau_i$, denote a  finite sequence of nodes generated from the distinct points
$z^{(k)}_\ell$, $1\leq \ell\leq n$, each of them repeated $\tau_\ell$ times.  Moreover, for a given $\eta>0$ let
$\gamma_\eta$ be the lemniscata defined by
$\gamma_\eta=\{z\in \mathbb C \colon |\prod_{\ell=1}^n(z-z^{(k)}_\ell)^{\tau_\ell}|=\eta^N\}$.
It is worth noticing that
\[
\prod_{\ell=1}^{n}(z-z^{(k)}_\ell)^{\tau_\ell}=\xi_k
(-z)^{\lceil \frac N2 \rceil}
B_{N,\Gamma}^{(k)}(-z),
\]
for a suitable $\xi_k\in \mathbb C$. Let us introduce the infinite sequence
$\{\tilde \beta_i\}_{i\in \mathbb N}$ obtained by cyclically repeating each $\beta_i$, i.e.,
$\tilde \beta_i=\beta_{{\rm mod}(i,N)}$, $i\geq 1$. For $m_1\leq m_2$  let  $f[\{\tilde \beta_i\}_{i=m_1}^{m_2}]$
be the divided difference of the meromorphic   function
$f(z) :=z^{\ell}a_{N,N,\Gamma}^{(k)}(z)/B^{(k)}_{N,\Gamma}(z)$  on the set of points
$\tilde \beta_i$, $m_1\leq i\leq m_2$, where $\ell \in \mathbb Z$ is a given fixed integer.
Then the relation
\[
f(z)= \sum_{i=1}^{+\infty} f[\{\tilde \beta_s\}_{s=1}^{i}]\prod_{j=1}^{i-1}(z-\tilde \beta_j)
\]
holds in the following sense: the partial sums of the Newton series converge uniformly to
$f(z)$ in any closed set lying in the interior of $\gamma_\eta$ for any $\eta>0$ such that
$f$ is analytic in $\gamma_\eta$.
Since
\[
\sum_{i=1}^{+\infty} f[\{\tilde \beta_s\}_{s=1}^{i}]\prod_{j=1}^{i-1}(z-\tilde \beta_j)=
 \sum_{i=1}^{N} f[\{\tilde \beta_s\}_{s=1}^{i}]\prod_{j=1}^{i-1}(z-\tilde \beta_j) + B_{N,\Gamma}^{(k)}(-z) R_k(z),
\]
one deduces that
\[
a_{N,N,\Gamma}^{(k)}(z)= z^{-\ell}B^{(k)}_{N,\Gamma}(z)\sum_{i=1}^{N} f[\{\tilde \beta_s\}_{s=1}^{i}]\prod_{j=1}^{i-1}(z-\tilde \beta_j)
+z^{-\ell}B^{(k)}_{N,\Gamma}(z)B_{N,\Gamma}^{(k)}(-z) R_k(z),
\]
and, therefore,  in view of \eqref{prop:Egeneration}, we can set
\[
c_{N,\Gamma}^{(k)}(z)= z^{-\ell}\sum_{i=1}^{N} f[\{\tilde \beta_s\}_{s=1}^{i}]\prod_{j=1}^{i-1}(z-\tilde \beta_j),
\]
which is the shifted  Newton form of the Hermite interpolant of
$\frac{a_{N,N,\Gamma}^{(k)}(z)}{B^{(k)}_{N,\Gamma}(z)}$
 at the points $z_\ell^{(k)},\ \ell=1,\cdots,n$.

\smallskip \noindent
In this way for any value of $\ell$ we  may  determine a Laurent polynomial  satisfying \eqref{intgen}
 whose support lies in $[-\ell, -\ell+N-1]$.
If $N$ is odd,  then by choosing $\ell=\frac{N-1}{2}$ we obtain the
unique Laurent polynomial supported in $[\frac{1-N}{2}, \frac{N-1}{2}]$. By using  the symmetry of    both
the symbol and  the distribution of nodes  and from the uniqueness of the Laurent polynomial
we may conclude that this latter polynomial is symmetric.  The case where $N$ is even is a bit  more involved. In principle,
using the same arguments as above we find that the interpolating Laurent polynomial  is supported in
$[-\frac{N}{2}, \frac{N}{2}-1]$ or  $[-\frac{N}{2}+1, \frac{N}{2}]$.  However,  by expressing the polynomial in
the new variable $t:=z+z^{-1}$ we find that there exists  a uniquely determined symmetric
interpolating Laurent polynomial supported in
$[-\frac{N}{2}+1, \frac{N}{2}-1]$.  The precise statement is given below.

\begin{proposition}\label{sympoly}
The polynomial correction $c_{N,\Gamma}^{(k)}(z)$ is the unique symmetric Laurent
polynomial supported in $[-\lceil \frac{N}{2}\rceil+1 ,$ $\lceil \frac{N}{2}\rceil -1 ]$ which satisfies the conditions
\eqref{intgen} with $p, N$ such that $p=0$ if $N$ even and $p=-\frac12$ if $N$ odd.
\end{proposition}

\proof
We are looking for a polynomial of the form
\[
\displaystyle{c_{N,\Gamma}^{(k)}(z)= c_0 +\sum_{j=1}^{\lceil \frac{N}{2}\rceil -1}c_j (z^j+z^{-j})},
\]
which fulfills  the interpolation conditions \eqref{intgen}. From Proposition \ref{prop:oddsymm_EPR}
it follows that the condition in $z_\ell^{(k)}$ implies the same in $(z_\ell^{(k)})^{-1}$ and viceversa.
Hence, for $N$ even we  have to impose $N/2$ independent conditions, while for $N$ odd we have
$\tau_{\ell}=2j+1$ conditions at $z_\ell^{(k)}=1$ plus  $(N-1)/2-j$ independent conditions.  Since from Remark \ref{add1}
the conditions at  $z_\ell^{(k)}=1$ yield $j+1$ independent conditions, we  obtain $(N+1)/2=\lceil \frac{N}{2}\rceil$
conditions also for  the  odd case.
Defining $t:=z+z^{-1}$ let us introduce the functions $p_j(t)=z^j + z^{-j}$, $j\geq 0$. Such functions
are monic Chebyshev-like polynomials of degree $j$ which, starting from $p_0(t)=2, \, p_1(t)=t$, satisfy
the three-term recurrence relation
\[
p_{j+1}(t)=t p_j(t) -p_{j-1}(t), \qquad j\geq 1.
\]
Hence, by writing
\[
c_{N,\Gamma}^{(k)}(z)=\frac{c_0}{2} p_0(t) +\sum_{j=1}^{\lceil \frac{N}{2}\rceil -1}c_j p_j(t)=\psi(t),
\]
the proof follows from the existence and uniqueness of the interpolating polynomial $\psi(t)$ on  the  considered
set of nodes.
\qed

\bigskip \noindent
The results given so far immediately generalize  to the case where we consider a subset ${\tilde \Gamma}$ of $\Gamma$.

\begin{theorem}\label{Theo:menocond} Let $\Gamma$  and
${\tilde \Gamma} \subset \Gamma$ be symmetric sets of cardinality $N$ and $M$, respectively, with $M$ and $N$ of the same parity.
 Moreover let $EP_{\Gamma}$, $EP_{\tilde \Gamma}$ be the corresponding
sets of exponential polynomials, and assume $p=0$ in case $N$ and $M$ are both even, $p=-\frac12$ in case $N$ and $M$ are both odd.
Then there exists a unique symmetric Laurent polynomial $c_{M,\Gamma}^{(k)}(z)$ supported in $[-\lceil \frac{M}{2}\rceil+1 ,
\lceil \frac{M}{2}\rceil -1 ]$ satisfying for $z^{(k)}_\ell:=e^{\frac{-\theta_\ell}{2^{k+1}}}$, $(\theta_\ell,\tau_\ell)\in {\tilde \Gamma}$, $\ell=1,...,m$,
the generalized interpolation conditions
\[
\frac{d^s\, c_{M,\Gamma}^{(k)}(z^{(k)}_\ell)}{dz^{s}}=
\sum_{i=0}^s \binom{s}{i}  v_{\ell, i}\frac{d^{s-i}\, G^{(k)}_{N,\Gamma}(z^{(k)}_\ell)}{dz^{s-i}}, \quad
 \begin{array}{l}
      \ell=1,\dots,m, \quad m\leq n,\\
      s=0,\dots,\tau_\ell-1,
\end{array}
\]
where $v_{\ell, i}:=2\left(z^{(k)}_\ell\right)^{p-i}\,\displaystyle{\prod_{r=0}^{i-1}(p-r)}$ and $\ G^{(k)}_{N,\Gamma}(z):=\frac 1{B^{(k)}_{N,\Gamma}(z)}$.

\end{theorem}

\medskip \noindent
For the effective construction of the polynomial $c_{M,\Gamma}^{(k)}(z)$ we can proceed as follows.
By setting
\begin{equation}\label{eq:cN_psi}
c_{M,\Gamma}^{(k)}(z)=\psi(z+z^{-1}),
\end{equation}
using the Fa\`a di Bruno's Formula we get that, for $s=1,\dots,\tau_\ell-1$,
\begin{equation}\label{eq:cN_psi_der}
\frac{d^s\, c_{M,\Gamma}^{(k)}(z^{(k)}_\ell)}{dz^{s}}=
\sum_{r=1}^s\frac{d^r\, \psi \left(z^{(k)}_\ell + (z^{(k)}_\ell)^{-1} \right)}{dz^{r}} \, A_{r,s}(z^{(k)}_\ell).
\end{equation}
Since
\[
A_{s,s}(z^{(k)}_\ell) =\left(1-\frac{1}{({z^{(k)}_\ell)}^2}\right)^s,
\]
we obtain that for $z^{(k)}_\ell\neq \pm 1$
the triangular system is invertible and, therefore, it enables the computation of the highest derivative of
$\psi(z+z^{-1})$ to be performed.  In this way $\psi(z+z^{-1})$  and  then $c_{M,\Gamma}^{(k)}(z)$ can be computed by using
the   customary Hermite interpolation formula.  For the case $z^{(k)}_\ell = 1$ it can be shown that
only the derivatives of $c_{M,\Gamma}^{(k)}(z)$ of even order  give  information about the derivatives of
$\psi(z+z^{-1})$. The case $z^{(k)}_\ell = -1$ cannot occur due to the definition of $z^{(k)}_\ell$.

\begin{proposition}\label{add2}
In the case $p=0$ and $M=N$ even, the subdivision symbol of the symmetric exponential pseudo-spline
$a_{N,N,\Gamma}^{(k)}(z)$ derived from Theorem \ref{Theo:menocond} is always interpolatory.
\end{proposition}

\proof
The proof exploits the fact that for $p=0$ and $N$ even from a subdivision symbol satisfying
\eqref{ERcond2} we are able to construct a subdivision symbol satisfying \eqref{ERcond1}, and viceversa.
In particular, in the case $p=0$ and $N$ even, from  conditions \eqref{intgen}
we find that $ a_{N,N,\Gamma}^{(k)}(z) - 2=  B^{(k)}_{N,\Gamma}(-z)\, q^{(k)}(z)$ for a certain  Laurent polynomial $q^{(k)}(z)$.
Hence,
\[
 a_{N,N,\Gamma}^{(k)}(z) = 2+  B^{(k)}_{N,\Gamma}(-z)\, q^{(k)}(z)
=B^{(k)}_{N,\Gamma}(z)\, c_{N,\Gamma}^{(k)}(z).
\]
Since the relation holds for any $z$ we also have
\[
2+  B^{(k)}_{N,\Gamma}(z)\, q^{(k)}(-z)
=B^{(k)}_{N,\Gamma}(-z)\, c_{N,\Gamma}^{(k)}(-z),
\]
which gives
\[
q^{(k)}(z) =-c_{N,\Gamma}^{(k)}(-z)
\]
and, hence,
\[
a_{N,N,\Gamma}^{(k)}(z)+a_{N,N,\Gamma}^{(k)}(-z)=2.
\]
\qed

\section{Convergence and regularity of exponential pseudo-spline subdivision schemes}\label{convergence}

The aim of this section is two-fold. First, we show that the strategy proposed in the previous section allows us to construct polynomial pseudo-splines in the stationary case. Second, we prove convergence and regularity of exponential pseudo-spline subdivision schemes.

\smallskip \noindent
To achieve the first goal we introduce a new set of $z$-functions, that are meant to be shifted B-spline symbols, of the form
$$
\bar{B}_N(z):=\frac{(z+1)^N}{2^{N-1}} \, z^{-\frac{N}{2}}=\left \{
\begin{array}{lll}
\frac{(z+z^{-1}+2)^{\rho}}{2^{2\rho-1}}, & \hbox{if} & N=2\rho\\
\frac{(z+z^{-1}+2)^{\rho+\frac12}}{2^{2\rho}}, & \hbox{if} & N=2\rho+1.
\end{array}
\right.
$$
Obviously, in case $N$ is even $\bar{B}_N(z)=B_N(z)$, whereas for $N$ odd $\bar{B}_N(z)=z^{\frac12}B_N(z)$.
We also need the following result whose proof is obtained by induction, following the lines of the proof of Proposition \ref{evensym_yoon}.

\begin{lemma}\label{ellam:yoon2}
The even-symmetric symbol $a(z)$ satisfies
$$\frac{d^r\, a(1)}{dz^r}=2\,\displaystyle{\prod_{i=0}^{r-1}\left(-\frac{1}{2}-i\right)},\quad r \geq 0,$$
if and only if  the associated odd-symmetric $z$-function $b(z)=z^{\frac12} \, a(z)$ satisfies
$$
\frac{d^r b(1)}{dz^r}=2 \delta_{r,0}, \qquad  r \geq 0.
$$
\end{lemma}

\smallskip \noindent By means of these preliminary results, we can prove the following proposition.

\begin{proposition}\label{prop:correct_stat}
For the order-$N$ B-spline symbol $B_{N}(z)=2^{-(N-1)} z^{-\lceil \frac N2 \rceil} \left(z+1\right)^N$, let $c_M(z)$ be the polynomial correction  such that
$$
 \sum_{i=0}^s \binom{s}{i}\frac{d^i\, c_{M}(1)}{dz^{i}} \frac{d^{s-i}\, B_{N}(1)}{dz^{s-i}}=
  2 \prod_{i=0}^{s-1}(p-i),\quad s=0,\dots,M-1,
$$
with $p=0$ if $N=2\rho$ and $p=-\frac12$ if $N=2\rho+1$. Then, ${a}_{M,N}(z)={B}_{N}(z) \, c_M(z)$ is the subdivision symbol of the polynomial pseudo-spline given in \cite[Section 6]{DHSS08}, that is
$$
{a}_{M,N}(z)=\left \{
\begin{array}{ll}
2 \,  \sigma^{\rho}(z)   \, \displaystyle{ \sum_{s=0}^{\left \lfloor  \frac{M-1}{2} \right \rfloor} \binom{\rho+s-1}{s} \delta^s(z)}, &  \hbox{if} \ N=2\rho,\smallskip \\
\frac{z+1}{z} \,  \sigma^{\rho}(z) \, \displaystyle{ \sum_{s=0}^{\left \lfloor  \frac{M-1}{2} \right \rfloor} \binom{\rho-\frac12+s}{s} \delta^s(z)}, &  \hbox{if} \ N=2\rho+1,
\end{array}
\right.
$$
where $\delta(z):=-\frac{(1-z)^2}{4z}$ and $\sigma(z):=\frac{(1+z)^2}{4z}$.
\end{proposition}

\proof
We start by observing that, since $c_{M}(z)=\frac{{a}_{M,N}(z)}{{B}_{N}(z)}=\frac{z^{\frac12}{a}_{M,N}(z)}{z^{\frac12}{B}_{N}(z)}$, in view of Lemma \ref{ellam:yoon2}, the polynomial correction $c_M(z)$ can be equivalently obtained by solving the linear system
$$
 \sum_{i=0}^s \binom{s}{i}\frac{d^i\, c_{M}(1)}{dz^{i}} \frac{d^{s-i}\, \bar{B}_{N}(1)}{dz^{s-i}}=
  2 v_s,\quad s=0,\dots,M-1
$$
with $\bar{B}_N(z)=2^{-(N-1)}\, z^{-\frac{N}{2}}\, (z+1)^N$ and $v_s=2\delta_{s,0}$.\\
We continue by taking
$$\bar{G}_N(z)=\frac{1}{\bar{B}_N(z)}=2^{2\rho-2p-1} (z+z^{-1}+2)^{-\rho+p}, \quad p \in \left \{0,-\frac12 \right \}.$$
Thus, using the Fa\`a di Bruno's formula (see \cite{WJ02} or \cite{Mortini}) we can write
$$
\frac{d^r \bar{G}_N(z)}{dz^r} =2^{2\rho-2p-1} \, \sum_{j=1}^r (-1)^j \binom{\rho-p+j-1}{j} j! (z+z^{-1}+2)^{-\rho+p-j} A_{j,r}(z),
$$
where
$$
A_{j,r}(z)=\sum_{\bq \in \bM^j, \, \vert \bq \vert=r} \frac{r!}{\bq!} \frac{\prod_{i=1}^j \left( \delta_{q_i,1} +(-1)^{q_i} q_i! z^{-(q_i+1)} \right)}{\prod_{i=1}^r N(\bq,i)!}
$$
with
$$
\bM^j=\{\bq=(q_1, q_2, ...,q_j) \in \NN^j, \, q_1 \geq q_2 \geq ... \geq q_j \geq 1  \}, \qquad \vert \bq \vert=q_1+...+q_j
$$
and $N(\bq,i)$ denoting the number of times the integer $i \in \NN$ appears in the $j$-tuple $\bq \in \NN^j$. Evaluating at $z=1$ we obtain
$$
\frac{d^r \bar{G}_N(1)}{dz^r} = \sum_{j=1}^r (-1)^j \binom{\rho-p+j-1}{j} j! 2^{-2j-1} A_{j,r}(1)
$$
so that, recalling \eqref{intgen}, for $s=0, \ldots, M-1$ we find %\left \lfloor  \frac{M-1}{2} \right \rfloor$ we find
$$
\frac{d^s c_M(1)}{dz^s} =\sum_{i=0}^s \binom{s}{i} v_i \sum_{j=1}^{s-i} (-1)^j \binom{\rho-p+j-1}{j} j! 2^{-2j-1} A_{j,s-i}(1).
$$
Hence,
\begin{eqnarray}\label{eq:rel1}
c_M(1)&=&1  \smallskip \\
\frac{d^s c_M(1)}{dz^s} &=& \displaystyle{ \sum_{j=1}^{s} (-1)^j \binom{\rho-p+j-1}{j} j! 2^{-2j} A_{j,s}(1)}, \quad s=1, \ldots, M-1.% \left \lfloor  \frac{M-1}{2} \right \rfloor.
\nonumber
\end{eqnarray}
On the other hand, recalling \eqref{eq:cN_psi}-\eqref{eq:cN_psi_der} we can write
\begin{eqnarray*}
c_M(z)&=& \psi(z+z^{-1}),\\
\frac{d^s c_M(z)}{d z^s}&=&\sum_{j=1}^{s} \frac{d^j \psi(z+z^{-1})}{dz^j} A_{j,s}(z), \quad s=1, \ldots, M-1,
\end{eqnarray*}
so that when evaluating at $z=1$ we obtain
\begin{eqnarray}\label{eq:rel2}
c_M(1)&=& \psi(2),\\
\frac{d^s c_M(1)}{d z^s}&=&\sum_{j=1}^{s} \frac{d^j \psi(2)}{dz^j} A_{j,s}(1), \quad s=1, \ldots, M-1.
\nonumber
\end{eqnarray}
In this way, by comparison of \eqref{eq:rel1} with \eqref{eq:rel2}, from all even $s$ we can find the values attained by all derivatives of $\psi$ at $2$, \ie
\begin{eqnarray*}
\psi(2)&=& 1,\\
\frac{d^j \psi(2)}{d z^j}&=&(-1)^{j} \binom{\rho-p+j-1}{j} j! 2^{-2j}, \quad j=1, \ldots, \left \lfloor  \frac{M-1}{2} \right \rfloor,
\end{eqnarray*}
and, exploiting the customary Hermite interpolation formula we can thus get the analytic expression of $c_M(z)$
$$
c_M(z) \equiv \psi(z+z^{-1})= \sum_{j=0}^{\left \lfloor  \frac{M-1}{2} \right \rfloor} \frac{(z+z^{-1}-2)^j}{j!} \, \frac{d^j \psi(2)}{d z^j}=\sum_{j=0}^{\left \lfloor  \frac{M-1}{2} \right \rfloor} \binom{\rho-p+j-1}{j} \delta^j(z),
$$
with $\delta(z)=-\, \frac{(1-z)^{2}}{4z}$. Making distinction between $N$ even and $N$ odd, and rewriting $B_N(z)$ in terms of $\sigma(z)=\frac{(1+z)^2}{4z}$, the claim is proven.
\qed

\smallskip \noindent
Now, in order to study the asymptotical behaviour of $a_{M,N,\Gamma}^{(k)}(z)$ when $k$ approaches infinity, we introduce the following definition.

\begin{definition} \label{def:asympt_simil} The sequence of subdivision
masks $\{\ba^{(k)},\ k\ge 0\}$ and $\{\ba\}$
are called \emph{asymptotically similar} if
$$\lim_{k\rightarrow +\infty}\|\ba^{(k)}-\ba\|_{\infty}=\lim_{k\rightarrow +\infty} \max_{i \in supp(\ba)} \vert \ra_i^{(k)}-\ra_i \vert =0,$$
or, equivalently in terms of symbols, if $\lim_{k \rightarrow +\infty} a^{(k)}(z)=a(z)$.
\end{definition}

\noindent
The following proposition proves the asymptotical similarity between the polynomial corrections $c_{M,\Gamma}^{(k)}(z)$ and $c_{M}(z)$.

\begin{proposition}\label{lem_inf}
As $k \rightarrow +\infty$ the exponential B-spline $B^{(k)}_{N,\Gamma}(z)$ converges to the polynomial B-spline $B_N(z)$, and the polynomial correction $c_{M,\Gamma}^{(k)}(z)$ approaches the symbol
\begin{equation}\label{eq:cm}
c_{M}(z)= \left \{
\begin{array}{ll}
\displaystyle{ \sum_{s=0}^{\left \lfloor  \frac{M-1}{2} \right \rfloor} \binom{\rho+s-1}{s} \delta^s(z)}, & \hbox{if} \ p=0 \ (\ie \, N=2\rho),\smallskip \\
\displaystyle{ \sum_{s=0}^{\left \lfloor  \frac{M-1}{2} \right \rfloor} \binom{\rho-\frac12+s}{s} \delta^s(z)}, & \hbox{if} \ p=-\frac12 \ (\ie \, N=2\rho+1),
\end{array}
\right.
\end{equation}
with $\delta(z):=-\, \frac{(1-z)^{2}}{4z}$.
\end{proposition}
\proof
We start by the simple observation that $\lim_{k \rightarrow +\infty} B^{(k)}_{N,\Gamma}(z)=B_N(z)$. Next, we first consider the case $p=0$ and $N$ even. From equation \eqref{eq:conditions} we have that
the Hermite interpolant of
$f(z):=\frac{a_{M,N,\Gamma}^{(k)}(z)}{B^{(k)}_{N,\Gamma}(z)}$
 at the points $z_\ell^{(k)},\ \ell=1,\cdots,m$ coincides with the Hermite interpolant of
$\chi^{(k)}(z):=\frac{2}{B^{(k)}_{N,\Gamma}(z)}$ at the same points.  Hence,
denoting by $\{ \beta_i^{(k)}\}_{i=1}^M$, $M=\sum_{\ell=1}^m \tau_{\ell}$, a  finite sequence of nodes generated from the distinct points
$z^{(k)}_\ell$, $1\leq \ell\leq m$, each of them repeated $\tau_\ell$ times,
and by $\tilde \beta_i^{(k)}=\beta^{(k)}_{{\rm mod}(i,M)}$, $i\geq 1$, we get
\[
f[\{\tilde \beta^{(k)}_s\}_{s=1}^{i}]=
\chi^{(k)}[\{\tilde \beta_s^{(k)}\}_{s=1}^{i}]=\frac{1}{2\pi{\tt i}} \int_C\frac{\chi^{(k)}(z) {\tt d}z}{\prod_{s=1}^i(z-\tilde \beta_i^{(k)})},
\]
where $C$ is a simple closed curve in the complex plane enclosing a simply connected region which contains the points
$\{\tilde \beta_s^{(k)}\}_{s=1}^{i}$. Therefore,
\[
\lim_{k\rightarrow +\infty} \chi^{(k)}[\{\tilde \beta_s^{(k)}\}_{s=1}^{i}]=\frac{1}{(i-1)!}
\frac{\partial^{i-1}}{\partial z^{i-1}}\chi^{(+\infty)}(z)\Big |_{z=1},
\]
which means that  as $k$ goes to infinity $c_{M,\Gamma}^{(k)}(z)$ approaches a shifted
 Newton form of the Hermite interpolant of
$\chi^{(+\infty)}(z)=\frac{2^{-(N-2)}}{z^{-\lceil \frac N2 \rceil} \left(
z+1\right)^N}$.

The remaining case $p=-1/2$ and $N$ odd reduces to the previous  analysis  by observing that
 $c_{M,\Gamma}^{(k)}(z)$ can be computed  from
\[
 z a_{M,N,\Gamma}^{(k)}(z^2) = B^{(k)}_{N,\Gamma}(z^2)\,  z c_{M,\Gamma}^{(k)}(z^2),
\]
by imposing the  generalized interpolation conditions \eqref{eq:conditions}  at the points $(z^{(k)}_{\ell})^\frac12,\ (z^{(k)}_{\ell})^{-\frac12}$
for the even symmetric symbol $\bar{a}_{M,N,\Gamma}^{(k)}(z)= z a_{M,N,\Gamma}^{(k)}(z^2)$.
\qed

\bigskip \noindent Collecting all previous results we finally arrive at an important asymptotical result.

\begin{coro}\label{prop:equivalence}
The exponential pseudo-spline subdivision masks are asymptotically similar to the polynomial pseudo-spline subdivision masks, \ie,
\begin{equation}\label{eq:limite}
\lim_{k\rightarrow +\infty} a_{M,N,\Gamma}^{(k)}(z)=a_{M,N}(z)
\end{equation}
with $a_{M,N}(z)$ denoting the well-known polynomial pseudo-spline symbol.
\end{coro}

\smallskip \noindent
Before proceeding, we recall that in   \cite{DongShen07} and in \cite{DDH10}  the authors prove convergence and regularity of the subdivision schemes associated with the symbol $a_{M,N}(z)$ with $N$ even and odd, respectively.

\smallskip \noindent
In the following we continue analyzing the values attained by $\{a_{M,N,\Gamma}^{(k)}(z),\ k\ge 0 \}$ and its derivatives at $z=-1$ and, in turn, we study the convergence and regularity of the associated subdivision schemes.

\begin{proposition}\label{prop:app_sum_rulesandconvergence}
Let $M>1$. The subdivision symbols $\{a_{M,N,\Gamma}^{(k)}(z),\ k\ge 0\}$ are such that
\begin{equation}\label{eq:approx}
\begin{array}{lll}
|a_{M,N,\Gamma}^{(k)}(1)-2|&=&O(2^{-kN}),\smallskip \\
\displaystyle \left | \frac {d^s}{dz^s} \, a_{M,N,\Gamma}^{(k)}(-1)\, \right |
      &=&O(2^{-k(N-s)}), \; \; s=0,\ldots,N-1,
\end{array}
\qquad k\to +\infty.
\end{equation}
Moreover, the non-stationary subdivision scheme with symbols  $\{a_{M,N,\Gamma}^{(k)}(z),\ k\ge 0\}$ converges and has the same regularity as the stationary one with symbol $a_{M,N}(z)$.
\end{proposition}

\proof
The proof of \eqref{eq:approx} is based on the recent results proven in \cite[Theorem 10]{ContiRomaniYoon2014} and is a direct consequence of the exponential polynomial generation properties of $\{a_{M,N,\Gamma}^{(k)}(z),\ k\ge 0\}$ and the asymptotical similarity of $\{\ba_{M,N,\Gamma}^{(k)},\ k\ge 0\}$ to $\{\ba_{M,N}\}$, previously shown in Proposition \ref{eq:limite}. Then, for the convergence and regularity result we can rely on
\cite{CharinaContiGuglielmiProtasov2014}.
\qed

\section{An application example}\label{Example}
This section contains an interesting example of a family of purely non-stationary symmetric exponential pseudo-splines. As far as we know this is the first derivation of purely non-stationary exponential pseudo-spline symbols to appear. In fact, the recently published paper \cite{NR14} merely discusses the interpolatory subcase of a family of exponential pseudo-splines which reproduces function spaces spanned by an arbitrary number of polynomials and just a pair of exponential polynomials.\\
Let $\theta \in \RR^+ \cup \ri[0,\pi)$ and for all $k\ge 0$ define $v^{(k)}=\frac{1}{2}(e^{\frac{\theta}{2^{k+1}}}+e^{-\frac{\theta}{2^{k+1}}})$.
We consider the exponential B-spline with $k$-level symbol
$$
B_{N,\Gamma}^{(k)}(z)=\frac{(z+z^{-1}+2v^{(k)})^{\rho}}{2^{2\rho-1} (v^{(k)})^\rho },
$$
with  $\rho \in \NN,$ $ \Gamma=\{ (\theta,\rho), (-\theta, \rho)\},$  and $N=\sharp \Gamma=2\rho,$
which is obtained from the general formulation with $n=2$, $z_1^{(k)}=e^{\frac{-\theta}{2^{k+1}}}$, $z_2^{(k)}=e^{\frac{\theta}{2^{k+1}}}$ and $\tau_1=\tau_2=\rho$.
The subdivision scheme with symbol $B_{N,\Gamma}^{(k)}(z)$ generates the space of exponential polynomials
\begin{equation}\label{eq:spanEP}
\hbox{span}\{ e^{\theta x}, \, e^{-\theta x}, \, x e^{\theta x}, \, x e^{-\theta x}, \, \cdots, x^{\rho-1} e^{\theta x}, \, x^{\rho-1} e^{-\theta x} \}
\end{equation}
and reproduces $\hbox{span}\{ e^{\theta x}, \, e^{-\theta x} \}$ with respect to the parameter shift $p=0$.
In order to apply Theorem \ref{Theo:menocond} with $m=n$, namely with $M=N$ (the only possibility we have here), we define
$$
G_{N,\Gamma}^{(k)}(z)=2^{2\rho-1} (v^{(k)})^{\rho} (z+z^{-1}+2v^{(k)})^{-\rho}
$$
and, using the Fa\`a di Bruno's formula we compute
$$
\frac{d^r G_{N,\Gamma}^{(k)}(z)}{d z^r}=2^{2\rho-1} (v^{(k)})^{\rho} \, \sum_{j=1}^r (-1)^j \binom{\rho+j-1}{j} j! (z+z^{-1}+2v^{(k)})^{-\rho-j} \, A_{j,r}(z)
$$
where $A_{j,r}(z)$ is the same as the one appearing in the proof of Proposition \ref{prop:correct_stat}.
Hence, for $\ell=1,2$
$$
G_{N,\Gamma}^{(k)}(z_{\ell}^{(k)})=\frac{1}{2}
$$
and
$$
\frac{d^r G_{N,\Gamma}^{(k)}(z_{\ell}^{(k)})}{d z^r}= \sum_{j=1}^r (-1)^j \binom{\rho+j-1}{j} j! 2^{-2j-1} (v^{(k)})^{-j} \, A_{j,r}(z_{\ell}^{(k)}), \quad r=1, \ldots, \rho-1.
$$
Therefore, from \eqref{intgen} we can write for $\ell=1,2$
\begin{equation}\label{eq:0}
\frac{d^s c_{N,\Gamma}^{(k)}(z_{\ell}^{(k)})}{d z^s}=\sum_{i=0}^s \binom{s}{i} v_{\ell,i} \sum_{j=1}^{s-i} (-1)^j \binom{\rho+j-1}{j} j! 2^{-2j-1} (v^{(k)})^{-j} \, A_{j,s-i}(z_{\ell}^{(k)}),
\end{equation}
with $s=0,...,\rho-1.$
Now, taking into account that when $p=0$ then $v_{\ell,i}=2 \delta_{i,0}$, $i=0, \ldots, s$, equation \eqref{eq:0} can be rewritten for $\ell=1,2$ as
\begin{eqnarray}
c_{N,\Gamma}^{(k)}(z_{\ell}^{(k)})=1, \hspace{7cm} \label{eq:vers1a}\\
\begin{array}{ll}
\displaystyle\frac{d^s c_{N,\Gamma}^{(k)}(z_{\ell}^{(k)})}{d z^s}= \sum_{j=1}^{s} (-1)^j \binom{\rho+j-1}{j} j! 2^{-2j} (v^{(k)})^{-j} \, A_{j,s}(z_{\ell}^{(k)}),
\\
 s=1, \ldots, \rho-1. \end{array}\label{eq:vers1b}
\end{eqnarray}
On the other hand, from \eqref{eq:cN_psi}-\eqref{eq:cN_psi_der} we have that
\begin{eqnarray*}
c_{N,\Gamma}^{(k)}(z)&=& \psi(z+z^{-1})\\
\frac{d^s c_{N,\Gamma}^{(k)}(z)}{d z^s}&=&\sum_{j=1}^{s} \frac{d^j \psi(z+z^{-1})}{dz^j} A_{j,s}(z), \quad s=1,...,\rho-1,
\end{eqnarray*}
which, when evaluated at $z_{\ell}^{(k)}$, $\ell=1,2$, yield
\begin{eqnarray}
c_{N,\Gamma}^{(k)}(z_{\ell}^{(k)})= \psi(2v^{(k)}) \hspace{5.5cm} \label{eq:vers2a}\\
\begin{array}{ll}
\displaystyle\frac{d^s c_{N,\Gamma}^{(k)}(z_{\ell}^{(k)})}{d z^s}=\sum_{j=1}^{s} \frac{d^j \psi(2v^{(k)})}{d z^j} A_{j,s}(z_{\ell}^{(k)}), \quad s=1,...,\rho-1,\end{array} \label{eq:vers2b}
\end{eqnarray}
due to the fact that $z_{\ell}^{(k)}+(z_{\ell}^{(k)})^{-1}=2v^{(k)}$ for $\ell=1,2$.
At this point, comparing \eqref{eq:vers1a} with \eqref{eq:vers2a} and \eqref{eq:vers1b} with \eqref{eq:vers2b}, we respectively obtain
\begin{eqnarray*}
\psi(2v^{(k)})&=&1,\\
\frac{d^j \psi(2v^{(k)})}{d z^j}&=& (-1)^j \binom{\rho+j-1}{j} j! 2^{-2j} (v^{(k)})^{-j}, \quad j=1,...,\rho-1.
\end{eqnarray*}
Using the customary Hermite interpolation formula we can thus compute
\begin{equation}\label{eq:cnz_form}
\begin{array}{lll}
\psi(z+z^{-1})&=&\displaystyle\sum_{j=0}^{\rho-1} \frac{(z+z^{-1}-2v^{(k)})^j}{j!} \frac{d^j \psi(2v^{(k)})}{d z^j}\\
&=&\displaystyle\sum_{j=0}^{\rho-1} (z+z^{-1}-2v^{(k)})^j
(-1)^j \binom{\rho+j-1}{j} 2^{-2j} (v^{(k)})^{-j},
\end{array}
\end{equation}
which provides the expression of the polynomial correction $c_{N,\Gamma}^{(k)}(z)$.\\
The obtained symbol $c_{N,\Gamma}^{(k)}(z)$ is the one that allows us to define the symmetric interpolatory symbol
$a_{N,N,\Gamma}^{(k)}(z)=B_{N,\Gamma}^{(k)}(z) c_{N,\Gamma}^{(k)}(z)$
which reproduces the space of exponential polynomials \eqref{eq:spanEP} with respect to the parameter shift $p=0$.
We can thus refer to $c_{N,\Gamma}^{(k)}(z)$ as to the polynomial correction transforming the approximating scheme having symbol $B_{N,\Gamma}^{(k)}(z)$ into the  interpolatory scheme with the same generation properties (see \cite{BCR10,CGR09,CGR11}).

\begin{rem}
Since, as previously observed, $\lim_{k \rightarrow +\infty} B_{N,\Gamma}^{(k)}=B_{2\rho}(z)$ and $\lim_{k \rightarrow +\infty}$ $c_{N,\Gamma}^{(k)}(z)=c_{2\rho}(z)$, then the derived family of exponential pseudo-splines with $k$-level symbol $a_{N,N,\Gamma}^{(k)}(z)$ can be considered a non-stationary extension of the family of interpolatory $(2\rho)$-point Dubuc-Deslauriers schemes \cite{DD89}. A general construction for families of interpolatory schemes reproducing exponential polynomials was originally proposed in \cite{DLL03}. However, to the best of our knowledge, an algebraic approach for deriving the subdivision symbols of exponential pseudo-splines (including as a special case all non-stationary variants of Dubuc-Deslauriers schemes) was never investigated before.
\end{rem}

\smallskip \noindent
For instance, note that when $\rho=2$, equation \eqref{eq:cnz_form} yields $c_{N,\Gamma}^{(k)}(z)=-\frac{1}{2v^{(k)}}z+2-\frac{1}{2v^{(k)}}z^{-1}$
and the resulting exponential pseudo spline is an interpolatory 4-point scheme with $k$-level mask
\begin{equation}\label{def:mask4pt}
\left \{-\frac{1}{16 (v^{(k)})^3}, \, 0, \, \frac{3(4(v^{(k)})^2 - 1)}{16 (v^{(k)})^3}, \, 1, \, \frac{3(4(v^{(k)})^2 - 1)}{16 (v^{(k)})^3}, \, 0, \, -\frac{1}{16 (v^{(k)})^3} \right \},
\end{equation}
while, when $\rho=3$, from equation \eqref{eq:cnz_form} we obtain $c_{N,\Gamma}^{(k)}(z)=\frac{3}{8(v^{(k)})^2} z^2 - \frac{9}{4v^{(k)}} z + \frac{3+16(v^{(k)})^2}{4(v^{(k)})^2} - \frac{9}{4v^{(k)}} z^{-1} +\frac{3}{8(v^{(k)})^2} z^{-2}$ and thus the
resulting exponential pseudo spline is an interpolatory 6-point scheme with $k$-level mask
\begin{equation}\label{def:mask6pt}
\begin{array}{ll}
\displaystyle \Big \{ \frac{3}{256 (v^{(k)})^5}, 0, -\frac{5(8(v^{(k)})^2 - 3)}{256 (v^{(k)})^5}, 0,
\frac{15 (8 (v^{(k)})^4 - 4(v^{(k)})^2 + 1)}{128 (v^{(k)})^5}, 1,
\smallskip \\
\quad \displaystyle \frac{15 (8 (v^{(k)})^4 - 4(v^{(k)})^2 + 1)}{128 (v^{(k)})^5}, 0, -
\frac{5(8(v^{(k)})^2 - 3)}{256 (v^{(k)})^5}, 0, \frac{3}{256 (v^{(k)})^5}  \Big \}.
\end{array}
\end{equation}

\begin{figure}[ht!]
\centering
{\includegraphics[trim= 10mm 25mm 5mm 30mm, clip, width=6.0cm]{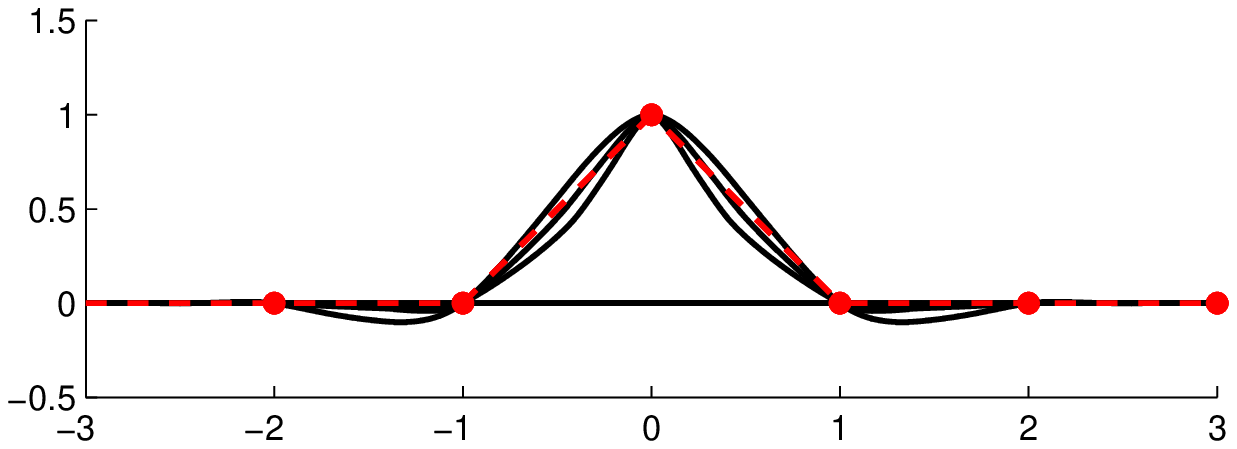}}\\
{\includegraphics[trim= 35mm 43mm 5mm 55mm, clip, width=13.5cm]{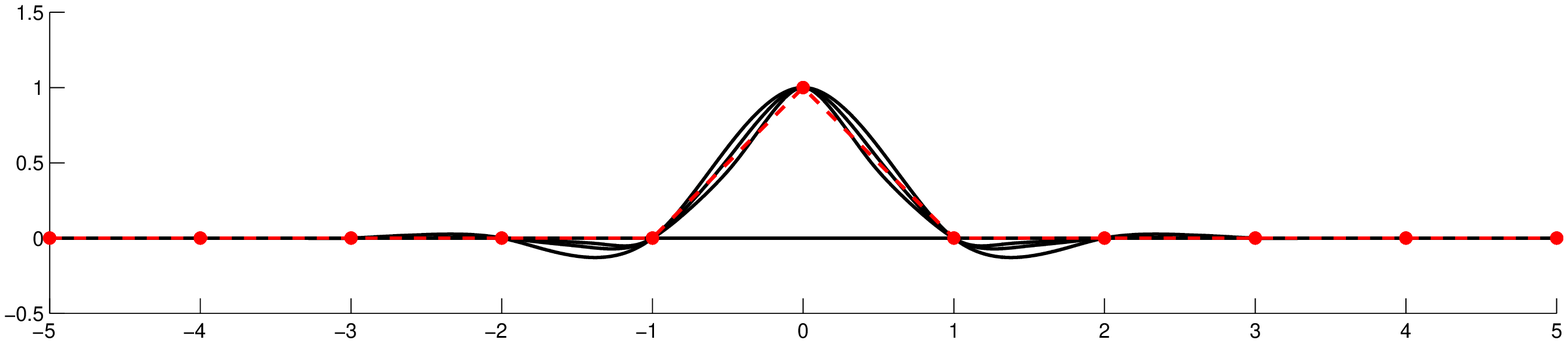}}
\caption{
Basic limit functions for the interpolatory 4-  and 6-point schemes with $k$-level
mask in \eqref{def:mask4pt} and \eqref{def:mask6pt} with $\theta=\ri$ (lower function), $\theta=\frac32$ (middle
function) and $\theta=2$ (upper function).}
\label{figB}
\end{figure}

\smallskip \noindent Figure \ref{figB} shows the graph of the basic limit functions for the interpolatory 4-point and 6-point schemes with $k$-level
mask in \eqref{def:mask4pt} and \eqref{def:mask6pt}, respectively.
Such interpolatory 4- and 6-point schemes are new non-stationary variants of the well-known Dubuc-Deslauriers schemes in \cite{DD89}. They differ from the ones previously proposed in \cite{BCR07,R09} for the space of exponential polynomials they reproduce.

\medskip

\section{Conclusions}\label{conclusion}
In this work we have proposed an algebraic strategy to derive the subdivision symbols of exponential pseudo-splines
from the subdivision symbols of exponential B-splines. The presented strategy is featured by the following key properties:
\begin{itemize}
\item it allows the user to pass from subdivision schemes generating a space of exponential polynomials to subdivision schemes reproducing the same space, or any desired of its subspaces;

\item it provides the subdivision symbols of minimal support that fulfill the set of conditions ensuring reproduction of the desired space of exponential polynomials;

\item it preserves the symmetry properties of the given exponential B-spline symbols;

\item it contains the stationary case of polynomial pseudo-splines as a special subcase.
\end{itemize}
Moreover, we have proved convergence and regularity of the non-stationary subdivision schemes obtained from the repeated application of exponential pseudo-spline symbols exploiting the property of asymptotical similarity to the stationary symbols of the well-known polynomial pseudo-splines.

\section*{Acknowledgements}
Lucia Romani acknowledges the support of MIUR-PRIN 2012 (grant 2012MTE38N).

\bibliographystyle{amsplain}

\end{document}